\documentclass[11pt]{article}
\usepackage{amsmath,amsfonts,amssymb}
\usepackage{theorem}
\usepackage[usenames]{color}

\theorembodyfont{\rmfamily}
\newtheorem{theorem}{Theorem}[section]
\newtheorem{proposition}[theorem]{Proposition}
\newtheorem{lemma}[theorem]{Lemma}
\newtheorem{definition}[theorem]{Definition}

\newtheorem{example}[theorem]{Example}
\newtheorem{remark}[theorem]{Remark}

\begin{document}
$\,$\vspace{10mm}

\begin{center}
{\textsf {\Huge A crystal theoretic method for finding}}
\vspace{2mm}\\
{\textsf {\Huge rigged configurations from paths}}
\vspace{13mm}\\
{\textsf {\LARGE Reiho Sakamoto}}
\vspace{2mm}\\
{\textsf {Department of Physics, Graduate School of Science,}}
\vspace{-1mm}\\
{\textsf {University of Tokyo, Hongo, Bunkyo-ku,}}
\vspace{-1mm}\\
{\textsf {Tokyo, 113-0033, Japan}}
\vspace{-1mm}\\
{\textsf {reiho@spin.phys.s.u-tokyo.ac.jp}}
\vspace{20mm}
\end{center}

\begin{abstract}
\noindent
The Kerov--Kirillov--Reshetikhin (KKR) bijection
gives one to one correspondences between
the set of highest paths and
the set of rigged configurations.
In this paper, we give a crystal theoretic
reformulation of the KKR map from
the paths to rigged configurations,
using the combinatorial $R$ and energy functions.
This formalism provides tool for analysis
of the periodic box-ball systems.
\end{abstract}
\pagebreak

\section{Introduction}
The Kerov--Kirillov--Reshetikhin (KKR) bijection \cite{KKR,KR,KSS}
gives the combinatorial one to one correspondences
between the set of highest weight elements of
tensor products of crystals \cite{Kas,KMN} (which we call highest paths)
and the set of combinatorial objects called the rigged configurations.
This bijection was originally introduced as essential tool
to derive new expression (called fermionic formulae)
 of the celebrated Kostka--Foulkas
polynomials \cite{Mac}.
Background of this expression is the
Bethe ansatz for the Heisenberg spin chain \cite{Bethe}
and, in this context, the rigged configurations form index set
of the eigenvalues and eigenvectors of the Hamiltonian.
To date the fermionic formulae have been extended to wider
class of representations, and proved in several cases
(see, e.g., \cite{O,Sch1,Lus,Nak,OS}
for the current status of the study).

Recently, the KKR bijection itself becomes active subject of the study.
Fundamental observation \cite{KOSTY} is that the KKR bijection is
inverse scattering transform of the box-ball systems, which is the
prototypical example of the ultradiscrete soliton systems
introduced by Takahashi--Satsuma \cite{TS,Tak}.
In this context, the rigged configurations are regarded as
action and angle variables for the box-ball systems.
This observation leads to derivation \cite{Sak,KSY} of general solutions
for the box-ball systems for the first time.

Therefore it is natural to ask what the representation
theoretic origin of the KKR bijection is.
Partial answer was given in the previous paper \cite{Sak},
and it is substantially used in the derivation in \cite{KSY}.
However, the formalism in \cite{Sak} works only for the
map from the rigged configuration to paths,
and extension of the formalism to the inverse direction
seems to have essential obstructions.
Up to now, crystal interpretation
for the map from the paths to rigged configurations
(which we denote by $\phi$) remains open.
Closely related problem is what the representation
theoretic origin of the mysterious combinatorial algorithm
of the original definition of $\phi$ is.
In fact, the formalism in \cite{Sak} gives alternative
representation theoretic map for $\phi^{-1}$
while it does not give meanings to the combinatorial
procedures like vacancy numbers or singular rows.
We remark that in Section 2.7 of \cite{KOSTY},
there is decomposition
of the $\mathfrak{sl}_n$ type $\phi$ into
successive computation of $\mathfrak{sl}_2$ type $\phi$.
However it finally uses combinatorial version of $\phi$,
hence it is not a complete crystal interpretation of $\phi$.

One of the main aims of the present article is to give a
crystal interpretation for the $\mathfrak{sl}_2$ type $\phi$
by clarifying the representation theoretic origin
of the original combinatorial procedure of $\phi$
(see Theorem \ref{s:main}).
In our formalism, the combinatorial procedure of $\phi$
is identified with differences of energy functions
called local energy distribution and indeed we can
read off all information about the rigged configurations
from them.
In terms of the box-ball systems,
the local energy distributions clarify which letters
of given path correspond to which soliton even if they
are in multiply scattering state.

Another aim of the paper is to provide a tool for
analysis of the periodic box-ball systems \cite{YT,YYT}.
In our $\mathfrak{sl}_2$ case formalism, there is remarkably
nice property (see Proposition \ref{s:no-crossing}).
Namely, the solitons appeared in the local energy distribution
are always separated from each other.
This leads to alternative version of our procedure
as given in Theorem \ref{s:thm:takagi}.
Importance of this reformulation is that when we apply
the formalism to find the action and angle variables of
the periodic box-ball systems \cite{KTT},
we do not need to cut paths and treat them as non-periodic paths.
This improves the inverse scattering formalism of \cite{KTT}
and theta function formulae of \cite{KS1,KS2}
in a sense that we treat paths genuinely as periodic.
As a byproduct, we give intuitive picture of the basic
operators which are the key to define angle variables
in \cite{KTT} (see Remark \ref{s:rem:pbbs}).
We remark that there is another approach
to the initial value problem of the periodic box-ball systems \cite{MIT}.
Although their combinatorial
method and our representation theoretic method are largely
different, it will be important to clarify
the relationship between these two approaches.

The present paper is organized as follows.
In Section \ref{s:sec:NY},
we review the combinatorial $R$ and energy functions
following \cite{NY}.
In Section \ref{s:sec:main}, we formulate our main results
(Theorem \ref{s:main} and
Theorem \ref{s:thm:takagi}).
In Section \ref{s:sec:proof}, we give proof of these theorems.
Section \ref{sec:summary} is summary.
In Appendix \ref{sec:kkr}, we recall the KKR bijection
and collect necessary facts, and in Appendix \ref{sec:T_l},
we collect necessary facts about the time evolution
operators $T_l$.

\section{Combinatorial $R$ and energy functions}\label{s:sec:NY}
In this section, we introduce necessary facts from
the crystal bases theory, namely, the combinatorial $R$
and energy functions.
Let $B_k$ be the crystal of $k$-fold symmetric powers of the
vector (or natural) representation of $U_q(\mathfrak{sl}_2)$.
As the set, it is
\begin{equation}
B_k=\{(x_1,x_2)\in\mathbb{Z}^2_{\geq 0}\,\vert\,
x_1+x_2=k\}.
\end{equation}
We usually identify elements of $B_k$ as the 
semi-standard Young tableaux
\begin{equation}
(x_1,x_2)=
\fbox{$\overbrace{1\cdots 1}^{x_1}\overbrace{2\cdots 2}^{x_2}$}\, ,
\end{equation}
i.e., the number of letters $i$ contained in a tableau is $x_i$.

For two crystals $B_k$ and $B_l$ of $U_q(\mathfrak{sl}_2)$,
one can define the tensor product
$B_k\otimes B_l=\{b\otimes b'\mid b\in B_k,b'\in B_l\}$.
Then we have a unique isomorphism $R:B_k\otimes B_l
\stackrel{\sim}{\rightarrow}B_l\otimes B_k$, i.e. a unique map
which commutes with actions of the Kashiwara operators
$\tilde{e}_i$, $\tilde{f}_i$.
We call this map combinatorial $R$
and usually write the map $R$ simply by $\simeq$.

In calculation of the combinatorial $R$,
it is convenient to use the diagrammatic technique
due to Nakayashiki--Yamada (Rule 3.11 of \cite{NY}).
Consider the two elements $x=(x_1,x_2)\in B_k$
and $y=(y_1,y_2)\in B_l$.
Then we draw the following diagram to express
the tensor product $x\otimes y$.
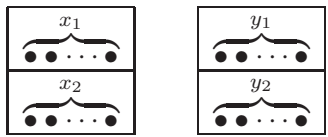
\begin{center}
\unitlength 12pt
\begin{picture}(10,4.5)
\multiput(0,0)(6,0){2}{
\multiput(0,0)(4,0){2}{\line(0,1){4}}
\multiput(0,0)(0,2){3}{\line(1,0){4}}
}
\put(0.5,0.2){$\overbrace{\bullet\bullet\cdots\bullet}^{x_{2}}$}
\put(0.5,2.2){$\overbrace{\bullet\bullet\cdots\bullet}^{x_1}$}
\put(6.5,0.2){$\overbrace{\bullet\bullet\cdots\bullet}^{y_{2}}$}
\put(6.5,2.2){$\overbrace{\bullet\bullet\cdots\bullet}^{y_1}$}
\end{picture}
\end{center}
The combinatorial $R$ matrix and energy function $H$ for
$x\otimes y\in B_k\otimes B_l$ (with $k\geq l$) are calculated by
the following rule.
\begin{enumerate}
\item
Pick any dot, say $\bullet_a$, in the right column and connect it
with a dot $\bullet_a'$ in the left column by a line.
The partner $\bullet_a'$ is chosen
{}from the dots
whose positions are higher than that of $\bullet_a$.
If there is no such a dot, we return to the bottom, and
the partner $\bullet_a'$ is chosen from the dots
in the lower row.
In the former case, we call such a pair ``unwinding,"
and, in the latter case, we call it ``winding."

\item
Repeat procedure (1) for the remaining unconnected dots
$(l-1)$ times.

\item
Action of the combinatorial $R$ matrix is obtained by
moving all unpaired dots in the left column to the right
horizontally.
We do not touch the paired dots during this move.

\item
The energy function $H$ is given by the number of unwinding pairs.
\end{enumerate}

The number of winding (or unwinding) pairs is sometimes called
the winding (or unwinding, respectively) number of tensor product.
It is known that the resulting combinatorial $R$ matrix
and the energy functions are not affected by the
order of making pairs
(\cite{NY}, Propositions 3.15 and 3.17).
In the above description, we only consider the case $k\geq l$.
The other case $k\leq l$ can be done by reversing the above
procedure, noticing the fact $R^2=\mathrm{id}$.
For more properties, including that the above
definition indeed satisfies the axiom, see \cite{NY}.

\begin{example}
Corresponding to the tensor product
$\fbox{1122}\otimes\fbox{122}$,
we draw the diagram given in the left hand side of:
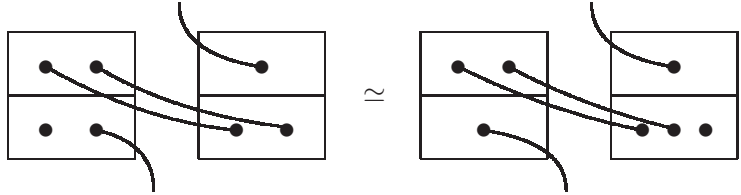
\begin{center}
\unitlength 12pt
\begin{picture}(24,6.5)(-0.2,-1)
\multiput(0,0)(6,0){2}{
\multiput(0,0)(4,0){2}{\line(0,1){4}}
\multiput(0,0)(0,2){3}{\line(1,0){4}}
}
\multiput(1.2,0.9)(1.6,0){2}{\circle*{0.4}}
\multiput(1.2,2.9)(1.6,0){2}{\circle*{0.4}}
\multiput(7.2,0.9)(1.6,0){2}{\circle*{0.4}}
\put(8.0,2.9){\circle*{0.4}}
\put(11.2,1.8){$\simeq$}
\thicklines
\qbezier(1.2,2.9)(4.2,1.2)(7.2,0.9)
\qbezier(2.8,2.9)(5.3,1.4)(8.8,1.0)
\qbezier(2.8,0.9)(4.6,0.5)(4.6,-1.0)
\qbezier(8.0,2.9)(5.4,3.3)(5.4,4.9)
\thinlines
\put(13,0){
\multiput(0,0)(6,0){2}{
\multiput(0,0)(4,0){2}{\line(0,1){4}}
\multiput(0,0)(0,2){3}{\line(1,0){4}}
}
\multiput(1.2,2.9)(1.6,0){2}{\circle*{0.4}}
\put(8.0,2.9){\circle*{0.4}}
\put(2.0,0.9){\circle*{0.4}}
\multiput(7.0,0.9)(1.0,0){3}{\circle*{0.4}}
\thicklines
\qbezier(1.2,2.9)(4.0,1.5)(7.0,0.9)
\qbezier(7.8,1.0)(5,1.7)(2.8,2.9)
\qbezier(8.0,2.9)(5.4,3.3)(5.4,4.9)
\qbezier(2.0,0.9)(4.6,0.5)(4.6,-1.0)
}
\end{picture}
\end{center}
By moving one unpaired dot to the right,
we obtain
\begin{equation}
\fbox{1122}\otimes\fbox{122}
\simeq
\fbox{112}\otimes\fbox{1222}\,.
\end{equation}
Since we have two unwinding pair, the energy function is
$H(\fbox{1122}\otimes\fbox{122})=2$.
\hfill$\square$
\end{example}

Consider the affinization of the crystal $B$.
As the set, it is
\begin{equation}
\mathrm{Aff}(B)=\{b[d]\, |\, b\in B,\, d\in\mathbb{Z}\}.
\end{equation}
Integers $d$ of $b[d]$ are called modes.
For the tensor product
$b_1[d_1]\otimes b_2[d_2]\in
\mathrm{Aff}(B_{k})\otimes\mathrm{Aff}(B_l)$,
we can lift the above definition of the
combinatorial $R$ as follows:
\begin{equation}
b_1[d_1]\otimes b_2[d_2]\stackrel{R}{\simeq}
b_2'[d_2-H(b_1\otimes b_2)]\otimes
b_1'[d_1+H(b_1\otimes b_2)],
\end{equation}
where $b_1\otimes b_2\simeq b_2'\otimes b_1'$
is the combinatorial $R$ defined in the above.
\begin{remark}
Piecewise linear formula to obtain the combinatorial $R$
and the energy function is also available \cite{HHIKTT}.
This is suitable for computer implementation.
For the affine combinatorial
$R:x[d]\otimes y[e]\simeq
\tilde{y}[e-H(x\otimes y)]\otimes
\tilde{x}[d+H(x\otimes y)]$,
we have
\begin{eqnarray}
&&\tilde{x}_i=x_i+Q_i(x,y)-Q_{i-1}(x,y),\qquad
\tilde{y}_i=y_i+Q_{i-1}(x,y)-Q_i(x,y),\nonumber\\
&&H(x\otimes y)=Q_0(x,y),\nonumber\\
&&Q_i(x,y)=\min (x_{i+1},y_i),
\end{eqnarray}
where we have expressed
$x=(x_1,x_2)$, $y=(y_1,y_2)$,
$\tilde{x}=(\tilde{x}_1,\tilde{x}_2)$ and
$\tilde{y}=(\tilde{y}_1,\tilde{y}_2)$.
All indices $i$ should be considered as
$i\in\mathbb{Z}/2\mathbb{Z}$.
\hfill$\square$
\end{remark}

\section{Local energy distribution and the KKR bijection}\label{s:sec:main}
In this section, we reformulate the combinatorial
procedure $\phi$ in terms of the energy functions
of crystal base theory.
See Appendix \ref{sec:kkr} for explanation of $\phi$.
In order to do this, it is convenient to express
actions of the combinatorial $R$ by vertex type diagrams.
First, we express the isomorphism of the combinatorial
$R$ matrix
\begin{equation}
a\otimes b_1\simeq b_1'\otimes a'
\end{equation}
and the corresponding value of the energy function
$e_1:=H(a\otimes b_1)$
by the following vertex diagram:
\begin{center}
\unitlength 12pt
\begin{picture}(4,4)
\put(0,2.0){\line(1,0){3.2}}
\put(1.6,1.0){\line(0,1){2}}
\put(-0.6,1.8){$a$}
\put(1.4,0){$b_1'$}
\put(1.4,3.2){$b_1$}
\put(3.4,1.8){$a'$}
\put(0.7,2.2){$e_1$}
\put(4.1,1.7){.}
\end{picture}
\end{center}
If we apply combinatorial $R$ successively as
\begin{equation}
a\otimes b_1\otimes b_2\simeq b_1'\otimes a'\otimes b_2
\simeq b_1'\otimes b_2'\otimes a'',
\end{equation}
with the energy function $e_2:=H(a'\otimes b_2)$,
then we express this by joining two vertices as follows:
\begin{center}
\unitlength 12pt
\begin{picture}(8,4)
\multiput(0,0)(4.2,0){2}{
\put(0,2.0){\line(1,0){3.2}}
\put(1.6,1.0){\line(0,1){2}}
}
\put(-0.6,1.8){$a$}
\put(1.4,0){$b_1'$}
\put(1.4,3.2){$b_1$}
\put(3.4,1.8){$a'$}
\put(5.6,3.2){$b_2$}
\put(5.6,0){$b_2'$}
\put(7.6,1.8){$a''$}
\put(0.7,2.2){$e_1$}
\put(4.9,2.2){$e_2$}
\put(8.3,1.7){.}
\end{picture}
\end{center}

\begin{definition}
For a given path $b=b_1\otimes b_2\otimes\cdots\otimes b_L$,
we define local energy $E_{l,j}$ by
$E_{l,j}:=H(u_l^{(j-1)}\otimes b_j)$.
Here, in the diagrammatic expression,
$u_l^{(j-1)}$ are defined as follows
(see also (\ref{s:def:T_l})
with convention $u_l^{(0)}:=u_l$).
\begin{center}
\unitlength 12pt
\begin{picture}(22,5)(0,-0.5)
\multiput(0,0)(5.8,0){2}{
\put(0,2.0){\line(1,0){4}}
\put(2,0){\line(0,1){4}}
}
\put(-0.9,1.8){$u_l$}
\put(0.5,2.3){$E_{l,1}$}
\put(1.7,4.2){$b_1$}
\put(1.7,-0.8){$b_1'$}
\put(4.2,1.8){$u_l^{(1)}$}
\put(6.3,2.3){$E_{l,2}$}
\put(7.5,4.2){$b_2$}
\put(7.5,-0.8){$b_2'$}
\put(10.0,1.8){$u_l^{(2)}$}
\multiput(11.5,1.8)(0.3,0){10}{$\cdot$}
\put(14.7,1.8){$u_l^{(L-1)}$}
\put(17,0){
\put(0,2.0){\line(1,0){4}}
\put(2,0){\line(0,1){4}}
}
\put(17.4,2.3){$E_{l,L}$}
\put(18.7,4.2){$b_L$}
\put(18.7,-0.8){$b_L'$}
\put(21.2,1.8){$u_l^{(L)}$}
\end{picture}
\end{center}
Here, we denote
$T_l(b)=b_1'\otimes b_2'\otimes\cdots\otimes b_L'$.
We define $E_{0,j}=0$ for all $1\leq j\leq L$.
We also use the following notation:
\begin{equation}
E_l:=\sum_{j=1}^L E_{l,j}.
\end{equation}
In other words, $u_l[0]\otimes b\stackrel{R}{\simeq}T_l(b)
\otimes u_l^{(L)}[E_l]$, where we have omitted
modes for $b$ and $T_l(b)$.

\hfill$\square$
\end{definition}
For a given path $b=b_1\otimes b_2\otimes\cdots\otimes b_L$
($b_i\in B_{\lambda_i}$), we create a path
$b'=b\otimes \fbox{1}^{\,\otimes\Lambda}$,
where $\Lambda >\lambda_1+\cdots +\lambda_L$.
Then we always have $u_l^{(L+\Lambda)}=u_l$ for
arbitrary $l$ (Proposition \ref{s:prop:ist} (1)).
In such a circumstance,
it is known that the sum $E_l$ is conserved
quantities of the box-ball system; $E_l(T_k(b'))=E_l(b')$.
The proof is based on successive application
of the Yang-Baxter equation (see Theorem 3.2 of \cite{FOY}
and section 3.4 of \cite{HHIKTT}).
However, for our purpose, we need more detailed
information such as $E_{l,j}$.

\begin{lemma}\label{s:lem:E=0,1}
For a given path $b=b_1\otimes b_2\otimes\cdots\otimes b_L$,
we have $E_{l,j}-E_{l-1,j}=0$ or 1,
for all $l>0$ and for all $1\leq j\leq L$.
\end{lemma}
{\bf Proof.}
First we give a proof when $l=1$, i.e., we show
$E_{1,j}-E_{0,j}=0$ or 1.
In this case, we have $u_l,u_l^{(i)}\in B_1$ and
$E_{0,j}=0$.
Since $H(x\otimes y)=0$ or 1 for all $x\in B_1$ and
all $y\in B_k$, the proof follows.

Now we consider possible values for $E_{l,j}-E_{l-1,j}$.
In order to do this, we show that the difference between
tableaux representations of $u_{l}^{(j)}$ and $u_{l-1}^{(j)}$
is only one letter.
More precisely, we show that if $u_{l-1}^{(j)}=(x_1,x_2)$,
then $u_{l}^{(j)}=(x_1+1,x_2)$ or
$u_{l}^{(j)}=(x_1,x_2+1)$.
We show this claim by induction on $j$.
For $j=0$ case, it is true because
$u_{l-1}^{(0)}=u_{l-1}=(l-1,0)$ and $u_{l}^{(0)}=u_{l}=(l,0)$,
by definition.
Suppose that the above claim holds for all $j<k$
for some $k$.
In order to compare $u_{l-1}^{(k)}$ and $u_l^{(k)}$,
consider the isomorphism
$u_{l-1}^{(k-1)}\otimes b_{k}\simeq b_{l-1,k}'\otimes u_{l-1}^{(k)}$
and
$u_l^{(k-1)}\otimes b_{k}\simeq b_{l,k}'\otimes u_l^{(k)}$.
By assumption, the difference between $u_{l-1}^{(k-1)}$
and $u_{l}^{(k-1)}$ is one letter.
Recall that in calculating the combinatorial $R$,
order of making pairs can be chosen arbitrary.
Therefore, in $u_l^{(k-1)}\otimes b_{k}$,
first we can make all pairs that appear in
$u_{l-1}^{(k-1)}\otimes b_{k}$,
and next we make remaining one pair.
This means the difference of number of unwinding pairs,
i.e., $E_{l,k}-E_{l-1,k}$ is 0 or 1.
To make the induction proceeds, note that this fact
means the difference between $u_{l-1}^{(k)}$
and $u_{l}^{(k)}$ is also one letter.
\hfill$\square$

\vspace{3mm}

The following theorem gives crystal theoretic
reformulation of the KKR map $\phi$.
See Appendix \ref{sec:kkr} for explanation of
the unrestricted rigged configurations.

\begin{theorem}\label{s:main}
Let $b=b_1\otimes b_2\otimes\cdots\otimes b_L
\in B_{\lambda_1}\otimes B_{\lambda_2}\otimes\cdots\otimes
B_{\lambda_L}$ be an arbitrary path.
$b$ can be highest weight or non-highest weight.
Set $N=E_1(b)$.
We determine the pair of numbers
$(\mu_1,r_1)$, $(\mu_2,r_2)$, $\cdots$, $(\mu_N,r_N)$
by the following procedure from Step 1 to Step 4.
Then the resulting $(\lambda ,(\mu ,r))$ coincides with
the (unrestricted) rigged configuration $\phi (b)$.
\begin{enumerate}
\item
Draw a table containing $(E_{l,j}-E_{l-1,j}=0,1)$
at the position $(l,j)$, i.e., at the $l$ th row and
the $j$ th column.
We call this table local energy distribution.
\item
Starting from the rightmost 1 in the $l=1$ st row,
pick one 1 from each successive row.
The one in the $(l+1)$ th row must be weakly right of
the one selected in the $l$ th row.
If there is no such 1 in the $(l+1)$ th row,
the position of the lastly picked 1 is called $(\mu_1,j_1)$.
Change all selected 1 into 0.
\item
Repeat Step 2 for $(N-1)$ times to further determine
$(\mu_2,j_2)$, $\cdots$, $(\mu_N,j_N)$
thereby making all 1 into 0.
\item
Determine $r_1,\cdots,r_N$ by
\begin{equation}\label{s:eq:rigging}
r_k=\sum_{i=1}^{j_k-1}\min (\mu_k,\lambda_i)
+E_{\mu_k,j_k}-2\sum_{i=1}^{j_k}E_{\mu_k,i}.
\end{equation}
\end{enumerate}
\hfill$\square$
\end{theorem}

Proof of Theorem \ref{s:main} will be given
in the next section.
As we will see in Proposition \ref{s:no-crossing},
the groups obtained in the above theorem
have no crossing with each other.
Therefore, when we search 1 of the $(l+1)$ th row
in the above Step 2,
we have at most one candidate, i.e.,
we can uniquely determine such 1.

\begin{example}\label{s:ex:mainth}
For example of Theorem \ref{s:main},
we consider the following path
\begin{equation}
b=
\fbox{1111}\otimes\fbox{11}\otimes\fbox{2}\otimes
\fbox{1122}\otimes\fbox{1222}\otimes\fbox{1}\otimes
\fbox{2}\otimes\fbox{22}
\end{equation}
Corresponding to Step 1,
the local energy distribution is given by
the following table
($j$ stands for column coordinate of the table).
\begin{center}
\begin{tabular}{l|cccccccc}
                 &1111&11&2 &1122&1222&1&2&22\\\hline
$E_{1,j}-E_{0,j}$&0   &0 &1 &0   &1   &0&1&0   \\
$E_{2,j}-E_{1,j}$&0   &0 &0 &1   &0   &0&0&1   \\
$E_{3,j}-E_{2,j}$&0   &0 &0 &1   &0   &0&0&0   \\
$E_{4,j}-E_{3,j}$&0   &0 &0 &0   &1   &0&0&0   \\
$E_{5,j}-E_{4,j}$&0   &0 &0 &0   &1   &0&0&0   \\
$E_{6,j}-E_{5,j}$&0   &0 &0 &0   &0   &0&0&1   \\
$E_{7,j}-E_{6,j}$&0   &0 &0 &0   &0   &0&0&0   \\
\end{tabular}
\end{center}
Following Step 2 and Step 3, letters 1
contained in the above table are found
to be classified into 3 groups,
as indicated in the following table.
\begin{center}
\begin{tabular}{l|cccccccc}
                 &1111&11&2 &1122&1222      &1  &2 &22        \\\hline
$E_{1,j}-E_{0,j}$&    &  &3 &    &$2^{\ast}$&   &1 &          \\
$E_{2,j}-E_{1,j}$&    &  &  &3   &          &   &  &$1^{\ast}$\\
$E_{3,j}-E_{2,j}$&    &  &  &3   &          &   &  &          \\
$E_{4,j}-E_{3,j}$&    &  &  &    &3         &   &  &          \\
$E_{5,j}-E_{4,j}$&    &  &  &    &3         &   &  &          \\
$E_{6,j}-E_{5,j}$&    &  &  &    &          &   &  &$3^{\ast}$\\
$E_{7,j}-E_{6,j}$&    &  &  &    &          &   &  &          \\
\end{tabular}
\end{center}
{}From the above table, we see that
the cardinalities of groups 1, 2 and 3 are
2, 1 and 6, respectively.
Also, in the above table,
positions of $(\mu_1,j_1)$, $(\mu_2,j_2)$
and $(\mu_3,j_3)$ are indicated by asterisks.
Their explicit locations are $(\mu_1,j_1)=(2,8)$,
$(\mu_2,j_2)=(1,5)$
and $(\mu_3,j_3)=(6,8)$.

Now we evaluate riggings $r_i$ according to
equation (\ref{s:eq:rigging}).
\begin{eqnarray*}
r_1&=&\sum_{i=1}^{8-1} \min(2,\lambda_i)
      +E_{2,8}-2\sum_{i=1}^8E_{2,i}\\
   &=&(2+2+1+2+2+1+1)+1-2(0+0+1+1+1+0+1+1)\\
   &=&2,\\
r_2&=&\sum_{i=1}^{5-1} \min(1,\lambda_i)
      +E_{1,5}-2\sum_{i=1}^5E_{1,i}\\
   &=&(1+1+1+1)+1-2(0+0+1+0+1)\\
   &=&1,\\
r_3&=&\sum_{i=1}^{8-1} \min(6,\lambda_i)
      +E_{6,8}-2\sum_{i=1}^8E_{6,i}\\
   &=&(4+2+1+4+4+1+1)+2-2(0+0+1+2+3+0+1+2)\\
   &=&1.
\end{eqnarray*}
Therefore we obtain
$(\mu_1,r_1)=(2,2)$, $(\mu_2,r_2)=(1,1)$ and $(\mu_3,r_3)=(6,1)$,
which coincide with the following computation
according to the original definition of $\phi$.
In the following, we use Young diagrammatic expression for
the rigged configurations, and we put riggings and
vacancy numbers on the right and on the left
of the corresponding rows of the configuration,
respectively.
\begin{center}
\unitlength 10pt
\begin{picture}(30,1.5)(0,0)
\put(0,0.2){$\emptyset$}
\put(2,0.2){$\emptyset$}
\put(3.5,0.2){$\stackrel{1}{\longrightarrow}$}
\multiput(6,0)(1,0){2}{\line(0,1){1}}
\multiput(6,0)(0,1){2}{\line(1,0){1}}
\put(9,0.2){$\emptyset$}
\put(10.5,0.2){$\stackrel{1}{\longrightarrow}$}
\put(6.5,0.5){\circle{0.4}}
\multiput(13,0)(1,0){3}{\line(0,1){1}}
\multiput(13,0)(0,1){2}{\line(1,0){2}}
\put(17,0.2){$\emptyset$}
\put(18.5,0.2){$\stackrel{1}{\longrightarrow}$}
\put(14.5,0.5){\circle{0.4}}
\multiput(21,0)(1,0){4}{\line(0,1){1}}
\multiput(21,0)(0,1){2}{\line(1,0){3}}
\put(26,0.2){$\emptyset$}
\put(27.5,0.2){$\stackrel{1}{\longrightarrow}$}
\put(23.5,0.5){\circle{0.4}}
\end{picture}
\end{center}
\begin{center}
\unitlength 10pt
\begin{picture}(30,2.5)(0,0)
\multiput(0,1)(1,0){5}{\line(0,1){1}}
\multiput(0,1)(0,1){2}{\line(1,0){4}}
\put(6,1.2){$\emptyset$}
\put(7.5,1.2){$\stackrel{1}{\longrightarrow}$}
\put(3.5,1.5){\circle{0.4}}
\multiput(10,1)(1,0){5}{\line(0,1){1}}
\multiput(10,1)(0,1){2}{\line(1,0){4}}
\multiput(10,0)(1,0){2}{\line(0,1){1}}
\put(10,0){\line(1,0){1}}
\put(16,1.2){$\emptyset$}
\put(17.5,1.2){$\stackrel{1}{\longrightarrow}$}
\put(10.5,0.5){\circle{0.4}}
\multiput(20,1)(1,0){5}{\line(0,1){1}}
\multiput(20,1)(0,1){2}{\line(1,0){4}}
\multiput(20,0)(1,0){3}{\line(0,1){1}}
\put(20,0){\line(1,0){2}}
\put(26,1.2){$\emptyset$}
\put(27.5,1.2){$\stackrel{2}{\longrightarrow}$}
\put(21.5,0.5){\circle{0.4}}
\end{picture}
\end{center}
\begin{center}
\unitlength 10pt
\begin{picture}(30,5.5)(0,-2)
\put(0,0){\line(1,0){1}}
\put(0,1){\line(1,0){2}}
\multiput(0,0)(1,0){2}{\line(0,1){1}}
\multiput(0,1)(1,0){3}{\line(0,1){1}}
\multiput(0,2)(1,0){5}{\line(0,1){1}}
\multiput(0,2)(0,1){2}{\line(1,0){4}}
\multiput(6,2)(1,0){2}{\line(0,1){1}}
\multiput(6,2)(0,1){2}{\line(1,0){1}}
\put(5.2,2.1){1}
\put(7.3,2.1){1}
\put(0.5,0.5){\circle{0.4}}
\put(6.5,2.5){\circle{0.4}}
\put(8.5,2.2){$\stackrel{2}{\longrightarrow}$}
\multiput(11,-1)(1,0){2}{\line(0,1){4}}
\multiput(11,-1)(0,1){2}{\line(1,0){1}}
\multiput(11,1)(0,1){1}{\line(1,0){2}}
\multiput(11,2)(1,0){5}{\line(0,1){1}}
\multiput(11,2)(0,1){2}{\line(1,0){4}}
\put(13,1){\line(0,1){1}}
\put(11.5,-0.5){\circle{0.4}}
\multiput(17,2)(1,0){3}{\line(0,1){1}}
\multiput(17,2)(0,1){2}{\line(1,0){2}}
\put(18.5,2.5){\circle{0.4}}
\put(16.2,2.2){2}
\put(19.3,2.2){2}
\put(20.5,2.2){$\stackrel{2}{\longrightarrow}$}
\end{picture}
\end{center}
\begin{center}
\unitlength 10pt
\begin{picture}(30,3.5)
\multiput(0,0)(0,1){2}{\line(1,0){2}}
\multiput(0,0)(1,0){2}{\line(0,1){4}}
\put(2,0){\line(0,1){1}}
\put(0,2){\line(1,0){2}}
\put(2,2){\line(0,1){1}}
\multiput(0,3)(0,1){2}{\line(1,0){4}}
\multiput(0,3)(1,0){5}{\line(0,1){1}}
\put(1.5,0.5){\circle{0.4}}
\multiput(6,3)(0,1){2}{\line(1,0){3}}
\multiput(6,3)(1,0){4}{\line(0,1){1}}
\put(8.5,3.5){\circle{0.4}}
\put(5.2,3.1){2}
\put(9.3,3.1){2}
\put(10.5,3.2){$\stackrel{1}{\longrightarrow}$}
\put(13,0)
{
\multiput(0,0)(0,1){2}{\line(1,0){3}}
\multiput(0,0)(1,0){2}{\line(0,1){4}}
\put(2,0){\line(0,1){1}}
\put(3,0){\line(0,1){1}}
\put(0,2){\line(1,0){2}}
\put(2,2){\line(0,1){1}}
\multiput(0,3)(0,1){2}{\line(1,0){4}}
\multiput(0,3)(1,0){5}{\line(0,1){1}}
\put(2.5,0.5){\circle{0.4}}
\multiput(6,3)(0,1){2}{\line(1,0){3}}
\multiput(6,3)(1,0){4}{\line(0,1){1}}
\put(5.2,3.1){3}
\put(9.3,3.1){2}
\put(10.5,3.2){$\stackrel{1}{\longrightarrow}$}
}
\end{picture}
\end{center}
\begin{center}
\unitlength 10pt
\begin{picture}(30,5.5)
\put(0,1)
{
\multiput(0,0)(0,1){2}{\line(1,0){4}}
\multiput(0,0)(1,0){2}{\line(0,1){4}}
\put(2,0){\line(0,1){1}}
\put(3,0){\line(0,1){1}}
\put(4,0){\line(0,1){1}}
\put(0,2){\line(1,0){2}}
\put(2,2){\line(0,1){1}}
\multiput(0,3)(0,1){2}{\line(1,0){4}}
\multiput(0,3)(1,0){5}{\line(0,1){1}}
\put(3.5,0.5){\circle{0.4}}
\multiput(6,3)(0,1){2}{\line(1,0){3}}
\multiput(6,3)(1,0){4}{\line(0,1){1}}
\put(5.2,3.1){3}
\put(9.3,3.1){2}
\put(10.5,3.2){$\stackrel{2}{\longrightarrow}$}
}
\put(13,1)
{
\multiput(0,-1)(1,0){2}{\line(0,1){1}}
\put(0,-1){\line(1,0){1}}
\multiput(0,0)(0,1){2}{\line(1,0){4}}
\multiput(0,0)(1,0){2}{\line(0,1){4}}
\put(2,0){\line(0,1){1}}
\put(3,0){\line(0,1){1}}
\put(4,0){\line(0,1){1}}
\put(0,2){\line(1,0){2}}
\put(2,2){\line(0,1){1}}
\multiput(0,3)(0,1){2}{\line(1,0){4}}
\multiput(0,3)(1,0){5}{\line(0,1){1}}
\put(0.5,-0.5){\circle{0.4}}
\multiput(6,2)(1,0){2}{\line(0,1){1}}
\put(6,2){\line(1,0){1}}
\put(6.5,2.5){\circle{0.4}}
\multiput(6,3)(0,1){2}{\line(1,0){3}}
\multiput(6,3)(1,0){4}{\line(0,1){1}}
\put(5.2,3.1){2}
\put(5.2,2.1){1}
\put(9.3,3.1){2}
\put(7.3,2.1){1}
\put(10.5,3.2){$\stackrel{2}{\longrightarrow}$}
}
\end{picture}
\end{center}
\begin{center}
\unitlength 10pt
\begin{picture}(30,5.5)
\put(0,1)
{
\multiput(0,-1)(1,0){3}{\line(0,1){1}}
\put(0,-1){\line(1,0){2}}
\multiput(0,0)(0,1){2}{\line(1,0){4}}
\multiput(0,0)(1,0){2}{\line(0,1){4}}
\put(2,0){\line(0,1){1}}
\put(3,0){\line(0,1){1}}
\put(4,0){\line(0,1){1}}
\put(0,2){\line(1,0){2}}
\put(2,2){\line(0,1){1}}
\multiput(0,3)(0,1){2}{\line(1,0){4}}
\multiput(0,3)(1,0){5}{\line(0,1){1}}
\put(1.5,-0.5){\circle{0.4}}
\multiput(6,2)(1,0){2}{\line(0,1){1}}
\put(6,2){\line(1,0){1}}
\put(9.5,3.5){\circle{0.4}}
\multiput(6,3)(0,1){2}{\line(1,0){4}}
\multiput(6,3)(1,0){5}{\line(0,1){1}}
\put(5.2,3.1){3}
\put(5.2,2.1){1}
\put(10.3,3.1){3}
\put(7.3,2.1){1}
\put(11.5,3.2){$\stackrel{2}{\longrightarrow}$}
}
\put(14,1)
{
\multiput(0,-1)(1,0){4}{\line(0,1){1}}
\put(0,-1){\line(1,0){3}}
\multiput(0,0)(0,1){2}{\line(1,0){4}}
\multiput(0,0)(1,0){2}{\line(0,1){4}}
\put(2,0){\line(0,1){1}}
\put(3,0){\line(0,1){1}}
\put(4,0){\line(0,1){1}}
\put(0,2){\line(1,0){2}}
\put(2,2){\line(0,1){1}}
\multiput(0,3)(0,1){2}{\line(1,0){4}}
\multiput(0,3)(1,0){5}{\line(0,1){1}}
\put(2.5,-0.5){\circle{0.4}}
\multiput(6,2)(1,0){2}{\line(0,1){1}}
\put(6,2){\line(1,0){1}}
\put(10.5,3.5){\circle{0.4}}
\multiput(6,3)(0,1){2}{\line(1,0){5}}
\multiput(6,3)(1,0){6}{\line(0,1){1}}
\put(5.2,3.1){2}
\put(5.2,2.1){1}
\put(11.3,3.1){2}
\put(7.3,2.1){1}
\put(12.5,3.2){$\stackrel{1}{\longrightarrow}$}
}
\end{picture}
\end{center}
\begin{center}
\unitlength 10pt
\begin{picture}(30,6.5)
\put(0,2)
{
\multiput(0,-1)(1,0){5}{\line(0,1){1}}
\put(0,-1){\line(1,0){4}}
\multiput(0,0)(0,1){2}{\line(1,0){4}}
\multiput(0,0)(1,0){2}{\line(0,1){4}}
\put(2,0){\line(0,1){1}}
\put(3,0){\line(0,1){1}}
\put(4,0){\line(0,1){1}}
\put(0,2){\line(1,0){2}}
\put(2,2){\line(0,1){1}}
\multiput(0,3)(0,1){2}{\line(1,0){4}}
\multiput(0,3)(1,0){5}{\line(0,1){1}}
\put(3.5,-0.5){\circle{0.4}}
\multiput(6,2)(1,0){2}{\line(0,1){1}}
\put(6,2){\line(1,0){1}}
\multiput(6,3)(0,1){2}{\line(1,0){5}}
\multiput(6,3)(1,0){6}{\line(0,1){1}}
\put(5.2,3.1){3}
\put(5.2,2.1){1}
\put(11.3,3.1){2}
\put(7.3,2.1){1}
\put(12.5,3.2){$\stackrel{1}{\longrightarrow}$}
}
\put(15,2)
{
\multiput(0,-2)(1,0){2}{\line(0,1){1}}
\put(0,-2){\line(1,0){1}}
\multiput(0,-1)(1,0){5}{\line(0,1){1}}
\put(0,-1){\line(1,0){4}}
\multiput(0,0)(0,1){2}{\line(1,0){4}}
\multiput(0,0)(1,0){2}{\line(0,1){4}}
\put(2,0){\line(0,1){1}}
\put(3,0){\line(0,1){1}}
\put(4,0){\line(0,1){1}}
\put(0,2){\line(1,0){2}}
\put(2,2){\line(0,1){1}}
\multiput(0,3)(0,1){2}{\line(1,0){4}}
\multiput(0,3)(1,0){5}{\line(0,1){1}}
\put(0.5,-1.5){\circle{0.4}}
\multiput(6,2)(1,0){2}{\line(0,1){1}}
\put(6,2){\line(1,0){1}}
\multiput(6,3)(0,1){2}{\line(1,0){5}}
\multiput(6,3)(1,0){6}{\line(0,1){1}}
\put(5.2,3.1){4}
\put(5.2,2.1){2}
\put(11.3,3.1){2}
\put(7.3,2.1){1}
\put(12.5,3.2){$\stackrel{2}{\longrightarrow}$}
}
\end{picture}
\end{center}
\begin{center}
\unitlength 10pt
\begin{picture}(30,7.5)
\put(0,3)
{
\multiput(0,-3)(1,0){2}{\line(0,1){2}}
\put(0,-3){\line(1,0){1}}
\put(0,-2){\line(1,0){1}}
\multiput(0,-1)(1,0){5}{\line(0,1){1}}
\put(0,-1){\line(1,0){4}}
\multiput(0,0)(0,1){2}{\line(1,0){4}}
\multiput(0,0)(1,0){2}{\line(0,1){4}}
\put(2,0){\line(0,1){1}}
\put(3,0){\line(0,1){1}}
\put(4,0){\line(0,1){1}}
\put(0,2){\line(1,0){2}}
\put(2,2){\line(0,1){1}}
\multiput(0,3)(0,1){2}{\line(1,0){4}}
\multiput(0,3)(1,0){5}{\line(0,1){1}}
\put(0.5,-2.5){\circle{0.4}}
\multiput(6,1)(1,0){2}{\line(0,1){2}}
\put(6,2){\line(1,0){1}}
\put(6,1){\line(1,0){1}}
\put(6.5,1.5){\circle{0.4}}
\multiput(6,3)(0,1){2}{\line(1,0){5}}
\multiput(6,3)(1,0){6}{\line(0,1){1}}
\put(5.2,3.1){3}
\put(5.2,2.1){1}
\put(5.2,1.1){1}
\put(11.3,3.1){2}
\put(7.3,2.1){1}
\put(7.3,1.1){1}
\put(12.5,3.2){$\stackrel{2}{\longrightarrow}$}
}
\put(15,3)
{
\multiput(0,-4)(1,0){2}{\line(0,1){3}}
\put(0,-4){\line(1,0){1}}
\put(0,-3){\line(1,0){1}}
\put(0,-2){\line(1,0){1}}
\multiput(0,-1)(1,0){5}{\line(0,1){1}}
\put(0,-1){\line(1,0){4}}
\multiput(0,0)(0,1){2}{\line(1,0){4}}
\multiput(0,0)(1,0){2}{\line(0,1){4}}
\put(2,0){\line(0,1){1}}
\put(3,0){\line(0,1){1}}
\put(4,0){\line(0,1){1}}
\put(0,2){\line(1,0){2}}
\put(2,2){\line(0,1){1}}
\multiput(0,3)(0,1){2}{\line(1,0){4}}
\multiput(0,3)(1,0){5}{\line(0,1){1}}
\put(0.5,-3.5){\circle{0.4}}
\multiput(7,1)(0,1){2}{\line(1,0){1}}
\put(8,1){\line(0,1){1}}
\multiput(6,1)(1,0){2}{\line(0,1){2}}
\put(6,2){\line(1,0){1}}
\put(6,1){\line(1,0){1}}
\put(7.5,1.5){\circle{0.4}}
\multiput(6,3)(0,1){2}{\line(1,0){5}}
\multiput(6,3)(1,0){6}{\line(0,1){1}}
\put(5.2,3.1){2}
\put(5.2,2.1){2}
\put(5.2,1.1){2}
\put(11.3,3.1){2}
\put(7.3,2.15){1}
\put(8.3,1.1){2}
\put(12.5,3.2){$\stackrel{2}{\longrightarrow}$}
}
\end{picture}
\end{center}
\begin{center}
\unitlength 10pt
\begin{picture}(30,8.5)
\put(0,3)
{
\multiput(1,-4)(0,1){2}{\line(1,0){1}}
\put(2,-4){\line(0,1){1}}
\multiput(0,-4)(1,0){2}{\line(0,1){3}}
\put(0,-4){\line(1,0){1}}
\put(0,-3){\line(1,0){1}}
\put(0,-2){\line(1,0){1}}
\multiput(0,-1)(1,0){5}{\line(0,1){1}}
\put(0,-1){\line(1,0){4}}
\multiput(0,0)(0,1){2}{\line(1,0){4}}
\multiput(0,0)(1,0){2}{\line(0,1){4}}
\put(2,0){\line(0,1){1}}
\put(3,0){\line(0,1){1}}
\put(4,0){\line(0,1){1}}
\put(0,2){\line(1,0){2}}
\put(2,2){\line(0,1){1}}
\multiput(0,3)(0,1){2}{\line(1,0){4}}
\multiput(0,3)(1,0){5}{\line(0,1){1}}
\put(1.5,-3.5){\circle{0.4}}
\multiput(7,1)(0,1){2}{\line(1,0){1}}
\put(8,1){\line(0,1){1}}
\multiput(6,1)(1,0){2}{\line(0,1){2}}
\put(6,2){\line(1,0){1}}
\put(6,1){\line(1,0){1}}
\multiput(6,3)(0,1){2}{\line(1,0){6}}
\multiput(6,3)(1,0){7}{\line(0,1){1}}
\put(11.5,3.5){\circle{0.4}}
\put(5.2,3.1){1}
\put(5.2,2.1){2}
\put(5.2,1.1){3}
\put(12.3,3.1){1}
\put(7.3,2.15){1}
\put(8.3,1.1){2}
}
\end{picture}
\end{center}
In the above diagrams, newly added boxes are
indicated by circles ``$\circ$".
The reader should
compare the local energy distribution and
the above box adding procedure.
Then one will observe that the local energy distribution and
box addition on $\mu$ part have close relationships.
This relation will be established in Lemma
\ref{s:lem:e_distrib->conf}.
In other words, the original combinatorial procedure
for $\phi$ is embedded
into rather automatic applications of the combinatorial $R$
and energy functions.
\hfill$\square$
\end{example}

\begin{example}\label{s:ex:mainth2}
By using Theorem \ref{s:main}, we can easily grasp the
large scale structure of combinatorial
procedures of the KKR bijection
from calculations of the combinatorial $R$
and energy functions.
In order to show the typical example,
consider the following long path
(length 40).
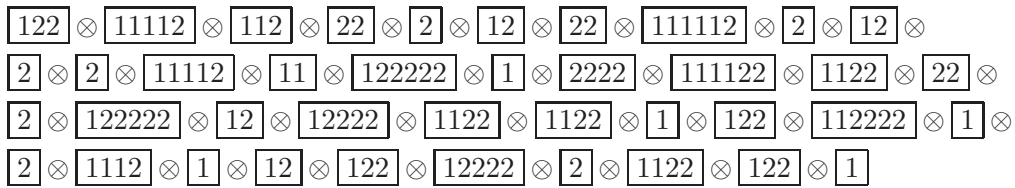
\begin{eqnarray*}
&&
\fbox{122}\otimes
\fbox{11112}\otimes
\fbox{112}\otimes
\fbox{22}\otimes
\fbox{2}\otimes
\fbox{12}\otimes
\fbox{22}\otimes
\fbox{111112}\otimes
\fbox{2}\otimes
\fbox{12}\otimes
\\
&&
\fbox{2}\otimes
\fbox{2}\otimes
\fbox{11112}\otimes
\fbox{11}\otimes
\fbox{122222}\otimes
\fbox{1}\otimes
\fbox{2222}\otimes
\fbox{111122}\otimes
\fbox{1122}\otimes
\fbox{22}\otimes
\\
&&
\fbox{2}\otimes
\fbox{122222}\otimes
\fbox{12}\otimes
\fbox{12222}\otimes
\fbox{1122}\otimes
\fbox{1122}\otimes
\fbox{1}\otimes
\fbox{122}\otimes
\fbox{112222}\otimes
\fbox{1}\otimes
\\
&&
\fbox{2}\otimes
\fbox{1112}\otimes
\fbox{1}\otimes
\fbox{12}\otimes
\fbox{122}\otimes
\fbox{12222}\otimes
\fbox{2}\otimes
\fbox{1122}\otimes
\fbox{122}\otimes
\fbox{1}
\end{eqnarray*}
Then, the local energy distribution takes
the following form.
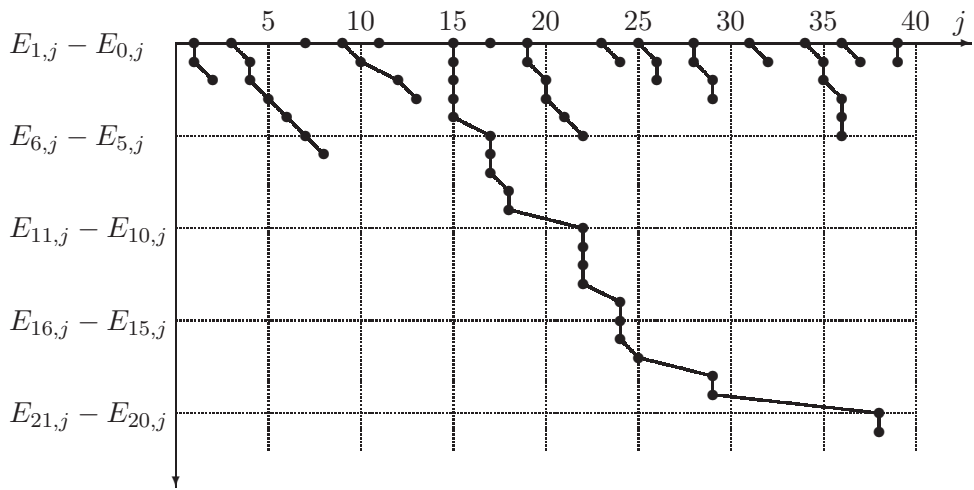
\begin{center}
\unitlength 7pt
\begin{picture}(45,25)(-6,-1.5)
\put(38,1){\circle*{0.5}}
\put(38,2){\circle*{0.5}}
\put(29,3){\circle*{0.5}}
\put(29,4){\circle*{0.5}}
\put(25,5){\circle*{0.5}}
\put(24,6){\circle*{0.5}}
\put(24,7){\circle*{0.5}}
\put(24,8){\circle*{0.5}}
\put(22,9){\circle*{0.5}}
\put(22,10){\circle*{0.5}}
\put(22,11){\circle*{0.5}}
\put(22,12){\circle*{0.5}}
\put(18,13){\circle*{0.5}}
\put(18,14){\circle*{0.5}}
\put(17,15){\circle*{0.5}}
\put(8,16){\circle*{0.5}}
\put(17,16){\circle*{0.5}}
\put(7,17){\circle*{0.5}}
\put(17,17){\circle*{0.5}}
\put(22,17){\circle*{0.5}}
\put(36,17){\circle*{0.5}}
\put(6,18){\circle*{0.5}}
\put(15,18){\circle*{0.5}}
\put(21,18){\circle*{0.5}}
\put(36,18){\circle*{0.5}}
\put(5,19){\circle*{0.5}}
\put(13,19){\circle*{0.5}}
\put(15,19){\circle*{0.5}}
\put(20,19){\circle*{0.5}}
\put(29,19){\circle*{0.5}}
\put(36,19){\circle*{0.5}}
\put(2,20){\circle*{0.5}}
\put(4,20){\circle*{0.5}}
\put(12,20){\circle*{0.5}}
\put(15,20){\circle*{0.5}}
\put(20,20){\circle*{0.5}}
\put(26,20){\circle*{0.5}}
\put(29,20){\circle*{0.5}}
\put(35,20){\circle*{0.5}}
\put(1,21){\circle*{0.5}}
\put(4,21){\circle*{0.5}}
\put(10,21){\circle*{0.5}}
\put(15,21){\circle*{0.5}}
\put(19,21){\circle*{0.5}}
\put(24,21){\circle*{0.5}}
\put(26,21){\circle*{0.5}}
\put(28,21){\circle*{0.5}}
\put(32,21){\circle*{0.5}}
\put(35,21){\circle*{0.5}}
\put(37,21){\circle*{0.5}}
\put(39,21){\circle*{0.5}}
\put(1,22){\circle*{0.5}}
\put(3,22){\circle*{0.5}}
\put(7,22){\circle*{0.5}}
\put(9,22){\circle*{0.5}}
\put(11,22){\circle*{0.5}}
\put(15,22){\circle*{0.5}}
\put(17,22){\circle*{0.5}}
\put(19,22){\circle*{0.5}}
\put(23,22){\circle*{0.5}}
\put(25,22){\circle*{0.5}}
\put(28,22){\circle*{0.5}}
\put(31,22){\circle*{0.5}}
\put(34,22){\circle*{0.5}}
\put(36,22){\circle*{0.5}}
\put(39,22){\circle*{0.5}}
\put(0,22){\vector(1,0){43}}
\put(0,22){\vector(0,-1){24}}
\multiput(5,0)(5,0){8}{
\multiput(0,0)(0,0.2){110}{\line(0,1){0.08}}
}
\multiput(0,2)(0,5){4}{
\multiput(0,0)(0.2,0){200}{\line(1,0){0.08}}
}
\put(4.6,22.7){5}
\put(9.2,22.7){10}
\put(14.2,22.7){15}
\put(19.2,22.7){20}
\put(24.2,22.7){25}
\put(29.2,22.7){30}
\put(34.2,22.7){35}
\put(39.2,22.7){40}
\put(42,22.7){$j$}
\put(-9,21.5){$E_{1,j}-E_{0,j}$}
\put(-9,16.5){$E_{6,j}-E_{5,j}$}
\put(-9,11.5){$E_{11,j}-E_{10,j}$}
\put(-9,6.5){$E_{16,j}-E_{15,j}$}
\put(-9,1.5){$E_{21,j}-E_{20,j}$}
\thicklines
\qbezier(39,22)(39,22)(39,21)
\qbezier(36,22)(36,22)(37,21)
\qbezier(34,22)(34,22)(35,21)
\qbezier(35,21)(35,21)(35,20)
\qbezier(35,20)(35,20)(36,19)
\qbezier(36,19)(36,19)(36,17)
\qbezier(31,22)(31,22)(32,21)
\qbezier(28,22)(28,22)(28,21)
\qbezier(28,21)(28,21)(29,20)
\qbezier(29,20)(29,20)(29,19)
\qbezier(25,22)(25,22)(26,21)
\qbezier(26,21)(26,21)(26,20)
\qbezier(23,22)(23,22)(24,21)
\qbezier(19,22)(19,22)(19,21)
\qbezier(19,21)(19,21)(20,20)
\qbezier(20,20)(20,20)(20,19)
\qbezier(20,19)(20,19)(22,17)
\qbezier(15,22)(15,22)(15,18)
\qbezier(15,18)(15,18)(17,17)
\qbezier(17,17)(17,17)(17,15)
\qbezier(17,15)(17,15)(18,14)
\qbezier(18,14)(18,14)(18,13)
\qbezier(18,13)(18,13)(22,12)
\qbezier(22,12)(22,12)(22,9)
\qbezier(22,9)(22,9)(24,8)
\qbezier(24,8)(24,8)(24,6)
\qbezier(24,6)(24,6)(25,5)
\qbezier(25,5)(25,5)(29,4)
\qbezier(29,4)(29,4)(29,3)
\qbezier(29,3)(29,3)(38,2)
\qbezier(38,2)(38,2)(38,1)
\qbezier(9,22)(9,22)(10,21)
\qbezier(10,21)(10,21)(12,20)
\qbezier(12,20)(12,20)(13,19)
\qbezier(3,22)(3,22)(4,21)
\qbezier(4,21)(4,21)(4,20)
\qbezier(4,20)(4,20)(8,16)
\qbezier(1,22)(1,22)(1,21)
\qbezier(1,21)(1,21)(2,20)
\end{picture}
\end{center}
In the above table, letters 1 in the local energy distribution
are represented by ``$\bullet$", and letters 0 are suppressed.
By doing Step 2 and Step 3, we obtain classifications
of letters 1.
In the above table, letters 1 belonging to the same
group are joined by thick lines.
We see there are 15 groups whose cardinalities are
3, 7, 1, 4, 1, 22, 1, 6, 2, 3, 4, 2, 6, 2, 2
from left to right, respectively.

By using equation (\ref{s:eq:rigging}), we obtain
the unrestricted rigged configuration as follows:
$(\mu_1,r_1)=(2,13)$, $(\mu_2,r_2)=(2,13)$,
$(\mu_3,r_3)=(6,15)$, $(\mu_4,r_4)=(2,12)$,
$(\mu_5,r_5)=(4,14)$, $(\mu_6,r_6)=(3,12)$,
$(\mu_7,r_7)=(2,11)$, $(\mu_8,r_8)=(6,5)$,
$(\mu_9,r_9)=(1,3)$, $(\mu_{10},r_{10})=(22,-17)$,
$(\mu_{11},r_{11})=(1,1)$, $(\mu_{12},r_{12})=(4,1)$,
$(\mu_{13},r_{13})=(1,1)$, $(\mu_{14},r_{14})=(7,-3)$
and $(\mu_{15},r_{15})=(3,-2)$.
The vacancy numbers for each row is
$p_{22}=-15$, $p_7=15$,
$p_6=19$, $p_4=21$,
$p_3=18$, $p_2=14$ and $p_1=10$.
Note that since the path in this example is not
highest weight, the resulting unrestricted rigged configuration
has negative values of the riggings and vacancy numbers.
\hfill$\square$
\end{example}

We have alternative form of Theorem \ref{s:main}.

\begin{theorem}\label{s:thm:takagi}
In the above Theorem \ref{s:main}, Step 2 can be replaced
by the following procedure (Step $2'$).
The resulting groups are the same as those obtained
in Theorem \ref{s:main} up to reordering in subscripts.
\begin{enumerate}
\item[$2'$.]
Pick one of the lowest 1 of the local energy distribution arbitrary,
and denote it by $(\mu_1,j_1)$.
Starting from $(\mu_1,j_1)$, choose one 1 from each row
successively as follows.
Assume that we have chosen 1 at $(l,k_{l})$.
Then $(l-1,k_{l-1})$ is the rightmost 1 among the part of
row $(l-1,1)$, $(l-1,2)$, $\cdots$, $(l-1,k_{l})$.
Change all selected 1 into 0.
\hfill$\square$
\end{enumerate}
\end{theorem}
Note that
comparing both Theorem \ref{s:main} and
Theorem \ref{s:thm:takagi},
the resulting $(\lambda ,(\mu ,r))$ can be different
in ordering of $\mu$.
However, this difference has no role in the KKR theory.
By the same reason, an ambiguity in choosing the lowest 1
in Step $2'$ brings no important difference.
Proof of Theorem \ref{s:thm:takagi}
will be given in the next section.

The formalism in Theorem \ref{s:thm:takagi}
is suitable for analysis of the periodic
box-ball system.
In particular, consider the case when
there are more than one longest group
in the local energy distribution.
Choose any successive longest groups and apply
the above Step $2'$ to these two groups.
Then due to non-crossing property of groups
(Proposition \ref{s:no-crossing}),
we can concentrate on the region between these two groups
and determine all groups between them
ignoring other part of the path.

\begin{remark}\label{s:rem:pbbs}
Let us remark how the above formalism works for
analysis of the periodic box-ball systems.
We concentrate on the path $b$ of the form $B_1^{\otimes L}$,
where number of \fbox{2} is equal to or less than
that of \fbox{1}.
We define $v_l\in B_l$
by the relation
$u_l\otimes b\simeq T_l(b) \otimes v_l$
with $T_l(b)\in B_1^{\otimes L}$.
Then we have
$v_l\otimes b\simeq \bar{T}_l(b)\otimes v_l$
with $\bar{T}_l(b)\in B_1^{\otimes L}$
(Proposition 2.1 of \cite{KTT}).
$\bar{T}_l$'s are the time evolution operator
of the periodic box-ball systems,
and $\bar{T}_1$ is simply the cyclic shift operator.

Consider the path $b^{\otimes N}=b\otimes\cdots\otimes b$.
Then, from the property of $v_l$, we have
$T_l(b^{\otimes N})=
T_l(b)\otimes \bar{T}_l(b)\otimes\cdots\otimes\bar{T}_l(b)$,
i.e., we can embed the periodic box-ball system
into the usual linear system with operator $T_l$.
Let us consider the local energy distribution
for $b^{\otimes N}$.
{}From the property
$v_l\otimes b\simeq \bar{T}_l(b)\otimes v_l$,
we see that under the right $(N-1)$ copies of $b$
in $b^{\otimes N}$, we have $(N-1)$ copies of
the same pattern of the local energy distribution.

Look at the local energy distribution below
the rightmost $b$.
In view of Theorem \ref{s:thm:takagi}
and the comments following it,
convenient way to find the structure of it
is as follows.
Instead of using $u_l$, we put $v_l$ on the left
of the path, and draw the local energy distribution.
In step $2'$ of Theorem \ref{s:thm:takagi},
we choose the rightmost 1 from $(l-1)$th row.
If there is not such 1, we return to the right end
of the $(l-1)$th row, and find such 1.

In such periodic extension of the local energy distribution,
we can always find a boundary of
successive two columns where none of groups
crossing the boundary.
By applying $\bar{T}_1$, we can move such boundary
to the left end of the path.
We assume that $b$ has already such property.
In our case, we can always do such procedure,
since by appropriate choice of $d$,
we can always make $\bar{T}_1^d(b)$
highest weight (such $d$ is not unique).
Then this $\bar{T}_1^d(b)$ meet the condition
(see Lemma C.1 of \cite{KTT}
and Lemma \ref{s:lem:e_distrib->conf}).

To summarize, by applying appropriate cyclic shifts,
we can always make $b^{\otimes N}$
whose local energy distribution is $N$ times repetition
of the pattern for single $b$.
On this property, we can apply the arguments of
\cite{KS1,KS2} (combined with the tau function of \cite{KSY})
to get the tau function in terms of
the ultradiscrete Riemann theta function.
More systematic treatment is given in \cite{KS3}.

Finally, we remark one thing without giving details
(see section 3.3 of \cite{KS3}).
Recall that there is ambiguity in the
choice of the cyclic shifts in the last paragraph.
Let $\bar{T}_1^d(b)$ and $\bar{T}_1^{d'}(b)$
are the two such possible choices
(we assume $d'=0$ for the sake of simplicity).
Consider the local energy distribution for $b$.
If the left $d$ columns contain a group of
cardinality $l$
(or, in other words, if the difference between
$b$ and $\bar{T}_1^d(b)$ is a soliton of length $l$),
then the riggings corresponding
to $b$ and $\bar{T}_1^d(b)$ differ by
the operator $\sigma_l$ called the slide
(see section 3.2 of \cite{KTT}
for definition of slides).
The slides $\sigma_l$ are closely related to the
period matrix of the tau functions of \cite{KS1,KS2}
(see section 4 of \cite{KTT}).
\hfill$\square$
\end{remark}

\section{Proof of Theorem \ref{s:main} and
Theorem \ref{s:thm:takagi}}\label{s:sec:proof}
For the proof of Theorem \ref{s:main},
we prepare some lemmas.
\begin{lemma}\label{s:lem:E=Q}
Let $(\lambda ,(\mu ,r))$ be the (unrestricted)
rigged configuration corresponding to the path
$b=b_1\otimes b_2\otimes\cdots\otimes b_L$.
Then we have (see (\ref{def:Q}) for
definition of $Q^{(1)}_l$)
\begin{equation}
E_l=Q^{(1)}_l.
\end{equation}
\end{lemma}
{\bf Proof.}
We consider the path
$b':=b\otimes\fbox{1}^{\,\otimes\Lambda}$,
where $\Lambda\gg |\lambda |$.
If $b'$ is not highest, apply Lemma \ref{s:lem:haiesutoka}
and we can use the same argument which is given below.
Since $E_l$ is conserved quantity on $b'$,
we have $E_l(T_\infty^{t_0}(b'))=E_l(b')$.
We take $t_0$ large enough with the
condition $\Lambda\geq t_0|\lambda |$
(the last inequality serves to assure that both
$T_\infty^{t_0}(b')$ and $b'$ contain the same
number of letters 2).
As we will see in the following,
$T_\infty^{t_0}(b')$ has simplified structure,
so that we can evaluate $E_l(T_\infty^{t_0}(b'))$
explicitly.

Now we use Proposition \ref{s:prop:ist}.
Since the actions of $T_\infty$ cause linear evolution
of riggings, we can assume
the (unrestricted) rigged configuration
corresponding to $T_\infty^{t_0}(b')$
as $(\lambda\cup (1^\Lambda ),(\mu, \bar{r}))$.
By the assumption $t_0\gg 1$,
these $\bar{r}$ have simple property.
Recall that in (\ref{s:prop:T_l:2}),
if we apply $T_\infty$ for one time, the rigging $r_i$
corresponding to the row $\mu_i$ becomes $r_i+\mu_i$.
Therefore the riggings $\bar{r}_i$ and $\bar{r}_j$
corresponding to the rows $\mu_i$ and $\mu_j$
satisfy $\bar{r}_i\gg \bar{r}_j$ if
$\mu_i>\mu_j$.

Using these observations, we determine the shape
of $T_\infty^{t_0}(b')$ from
$(\lambda\cup (1^\Lambda ),(\mu, \bar{r}))$.
By the assumption $t_0\gg 1$, all letters 2 in $T_\infty^{t_0}(b')$
are contained in $B_1^{\,\otimes\Lambda}$ part of the path.
Therefore, corresponding to the row $\mu_i$,
there is a soliton of the form $\fbox{2}^{\,\otimes\mu_i}$.
For example, in the following path,
\begin{equation}
\cdots\otimes
\fbox{1}\otimes\fbox{2}\otimes
\overbrace{
\fbox{1}\otimes\fbox{1}\otimes\cdots\otimes
\fbox{1}\otimes\fbox{1}
}^{\gg 1}
\otimes\fbox{2}\otimes\fbox{2}\otimes
\fbox{1}\otimes\fbox{1}\otimes\fbox{2}\otimes\fbox{2}\otimes
\fbox{1}\otimes\cdots
\end{equation}
there are one soliton of length 1 and two solitons
of length 2.
Since the riggings satisfy $\bar{r}_i\gg \bar{r}_j$ if
$\mu_i>\mu_j$,
the shorter solitons are located on the far left of
longer solitons (see the above example).

Assume there are solitons of the same length
such as $\fbox{2}^{\,\otimes\mu_2}\otimes\fbox{1}^{\,\otimes\sigma}
\otimes\fbox{2}^{\,\otimes\mu_1}$ ($\mu_1=\mu_2$).
Then we show $\sigma\geq\mu_1=\mu_2$.
Let the riggings corresponding to the rows $\mu_1$ and $\mu_2$
be $r_1$ and $r_2$, respectively.
In order to minimize $\sigma$, we choose $r_1=r_2$.
Now we consider $\phi^{-1}$ on rows $\mu_1$ and $\mu_2$.
Since we are assuming $t_0\gg 1$, we do not need to consider
the rows whose widths are different from $\mu_1$.
{}From $r_1=r_2$, rows $\mu_1$ and $\mu_2$ become simultaneously
singular, and we can choose one of them arbitrary.
We remove $\mu_1$ first.
While removing boxes from $\mu_1$ one by one,
the shortened row $\mu_1$ is always made singular,
and the rows whose lengths are shorter than $\mu_1$
are not singular.
Therefore we can remove entire row $\mu_1$ successively.
After removing row $\mu_1$, $Q^{(0)}_{\mu_2}$ decrease
by $\mu_1$ (note that the shape of the removed part of the
quantum space is $(1^{\mu_1})$),
and $Q^{(1)}_{\mu_2}$ is also decrease by $\mu_1$
(because of the removal of $\mu_1$).
Since the vacancy number is defined by
$Q^{(0)}_{\mu_2}-2Q^{(1)}_{\mu_2}$,
the vacancy number for the row $\mu_2$ increase by $\mu_1$
compared to the one calculated before removing $\mu_1$.
Therefore, in order to make the row $\mu_2$ singular again,
we have to remove extra $\mu_1$ boxes from the
quantum space, without removing boxes of $\mu$ part.
Hence we have $\sigma\geq\mu_1$, as requested.

Now we are ready to evaluate $E_l(T_\infty^{t_0}(b'))$.
{}From definition of the combinatorial $R$,
we have
\begin{eqnarray}
&&
\fbox{$\overbrace{11\cdots 1}^l$}\otimes
\fbox{2}^{\,\otimes m}\otimes
\fbox{1}^{\,\otimes M}\nonumber\\
&\simeq&
\fbox{1}^{\,\otimes\min (l,m)}\otimes
\fbox{2}^{\,\otimes\max (m-l,0)}\otimes
\fbox{$\overbrace{11\cdots 1}^{\max (l-m,0)}
\overbrace{22\cdots 2}^{\min (l,m)}$}\otimes
\fbox{1}^{\,\otimes M}\label{s:eq:combr}
\end{eqnarray}
As we have seen, if there is a soliton of length $m$,
there is always interval longer than $m$,
i.e., it has the form $\cdots\otimes\fbox{2}^{\,\otimes m}\otimes
\fbox{1}^{\,\otimes M}\otimes\cdots$ with $m\leq M$.
This makes \fbox{$11\cdots 122\cdots 2$} of (\ref{s:eq:combr})
into the form $\fbox{$11\cdots 1$}=u_l$
when it comes to the left of the next (or right) soliton.
Therefore, in order to evaluate $E_l(T_\infty^{t_0}(b'))$,
we only have to consider the situation
like (\ref{s:eq:combr}).
Noticing the fact that the energy function, i.e.,
unwinding number gains its value from the unwinding pair
$\fbox{$1\cdots$}\otimes\fbox{$2$}$
appearing in the each tensor product
(more precisely, $E_l$ gains $\min (l,m)$ corresponding to 
the procedure (\ref{s:eq:combr})),
the proof of lemma finishes.
\hfill$\square$

\vspace{3mm}

Combining the property of the combinatorial $R$ with
Lemma \ref{s:lem:E=Q},
the relationship between the local energy distribution
and the KKR bijection can be clarified as follows.

\begin{lemma}\label{s:lem:e_distrib->conf}
For the given path
$b=b_1\otimes\cdots\otimes b_i\otimes\cdots\otimes b_L$,
draw local energy distribution.
Within the $i$ th column, denote the locations of 1
as $(j_1,i)$, $(j_2,i)$, $\cdots$, $(j_k,i)$
($j_1<j_2<\cdots <j_k$).
Consider the calculation of
$\phi (b)=(\lambda,(\mu,r))$.
During the whole process of $\phi (b)$,
when we create part of the (unrestricted) rigged configuration
from $b_i$ of $b$, we add boxes to columns
$j_1,j_2,\cdots,j_k$ of $\mu$ in this order.
\end{lemma}
{\bf Proof.}
Denoting $b_i=(x_1,x_2)$, let us define
$b_{i,s}=(0,s)$ for $s\leq x_2$.
Consider the path
$b_s=b_1\otimes b_2\otimes\cdots\otimes b_{i,s}$,
and draw local energy distribution for this $b_{s}$.
In the local energy distribution,
from the 1 st column to $i-1$ th column are identical to
the ones in the case for the original $b$.
On the other hand, from $i+1$ th column to
$L$ th column, local energies are all 0.
These are obvious from the construction of $b_s$.

Now consider the $i$ th column of the local energy
distribution for $b_s$.
Then we show that the $i$ th column is obtained from the one
corresponding to the original $b$ by making $s$ letters 1
from the top as it is, and letting all other letters be 0.
This follows from the property of the combinatorial $R$,
that is, the order of making pairs of dots  do not affect
the resulting image of the combinatorial $R$.
In fact, in calculating $E_{l,i}$, we have
$u_l^{(i-1)}\otimes b_{i,s}$.
Compare this with $u_l^{(i-1)}\otimes b_{i}$.
Then we can make pair of dots in $u_l^{(i-1)}\otimes b_{i}$
such that first we join unwinding pairs, and then
we join winding pairs.
Note that all letters 1 contained in $b_i$ here
cannot contribute as unwinding pairs.
Therefore, we see that when we consider
$u_l^{(i-1)}\otimes b_{i,s}$,
$E_{l,i}$ ($l=1,2,\cdots$) are the same with
$u_l^{(i-1)}\otimes b_{i}$ case
up to the first $s$ unwinding pairs, and we do not
have the rest of the unwinding pairs.
This verifies the assertion for $i$ th column of the local energy
distribution for $b_s$.

Compare the local energy distribution for
$b_1\otimes\cdots\otimes b_{s-1}$ and
$b_1\otimes\cdots\otimes b_s$.
Then, from the above observation,
there is extra one 1 at $(j_s,i)$.
Now we apply the relation $E_l=Q^{(1)}_l$
(Lemma \ref{s:lem:E=Q}) to both
$b_1\otimes\cdots\otimes b_{s-1}$ and
$b_1\otimes\cdots\otimes b_s$.
Then we see that the letter 1 at $(j_s,i)$
corresponds to addition of one box at column $j_s$ of $\mu$
of $\phi (b_{s-1})$.
Since $\mu$ part of $\phi (b)$ is obtained by adding boxes
to $\mu$ recursively as
$\phi (b_{1}),\phi (b_2),\cdots,\phi (b_{i,1}),\phi (b_{i,2}),\cdots$,
this gives the proof of lemma.
\hfill$\square$

\begin{lemma}\label{s:junbi}
Let $b=b_1\otimes b_2\otimes\cdots\otimes b_L
\in B_{\lambda_1}\otimes B_{\lambda_2}\otimes\cdots\otimes
B_{\lambda_L}$ be an arbitrary path.
$b$ can be highest weight or non-highest weight.
Set $N=E_1(b)$.
We determine the numbers
$\mu_1$, $\mu_2$, $\cdots$, $\mu_N$
by the following procedure from Step 1 to Step 3.
\begin{enumerate}
\item
Draw a table containing $(E_{l,j}-E_{l-1,j}=0,1)$
at the position $(l,j)$, i.e., at the $l$ th row and
the $j$ th column.
We call this table local energy distribution.
\item
Starting from the rightmost 1 in the $l=1$ st row,
pick the nearest 1 from each successive row.
The one in the $(l+1)$ th row must be weakly right of
the one selected in the $l$ th row.
If there is no such 1 in the $(l+1)$ th row,
the position of the lastly picked 1 is called $(\mu_1,j_1)$.
Change all selected 1 into 0.
\item
Repeat Step 2 for $(N-1)$ times to further determine
$(\mu_2,j_2)$, $\cdots$, $(\mu_N,j_N)$
thereby making all 1 into 0.
Then $\mu$ coincides with $\mu$ of
the (unrestricted) rigged configuration
$\phi (b)=(\lambda ,(\mu ,r))$.
\end{enumerate}
\end{lemma}
{\bf Proof.}
We first interpret Step 2 in terms of the original combinatorial
procedure $\phi$.
In Step 2, we choose the rightmost 1 of the first row
of local energy distribution.
From Lemma \ref{s:lem:e_distrib->conf},
this 1 corresponds to the leftmost box of the lastly
created row of $\mu$.
Suppose we choose letters 1 up to $l$ th row
according to Step 2.
Next we choose 1 in $l+1$ th row,
whose position is weakly right of 1 in $l$ th row.
Since the lastly created row grows by adding boxes one by one
during the procedure $\phi$,
this means that these two 1 at $l$ th and $l+1$ th rows
of the local energy distribution
belong to the same row (lastly created row) of $\mu$.
Note that
if there are more than one row with the same length $l$,
we can always add a box to the lastly created row,
since it has maximal riggings among the rows
with the same length.
This follows from the fact that after adding a
box at $l$ th column of lastly created row of $\mu$,
the row is made to be singular, i.e., the row is assigned
the maximal possible rigging for the row with length $l$.
To summarize, Step 2 ensures us to identify all 1 in the
local energy distribution which correspond to
the lastly created row of $\mu$.

In Step 3, we do the same procedure for the rest
of 1 in the local energy distribution.
Since we omit all letters 1 which are already
identified with some rows of $\mu$,
we can always use Step 2 to determine the next row.
Therefore, Step 3 ensures us to identify all 1
in the local energy distribution with the rows of $\mu$.
\hfill$\square$

\vspace{3mm}

For the proof of Theorem \ref{s:main},
we show that the groups of letters 1 obtained in
Lemma \ref{s:junbi} have simplified structure.

\begin{proposition}\label{s:no-crossing}
The groups of letters 1 obtained in Step 2 and Step 3
of Lemma \ref{s:junbi} have no crossing with each other.
\end{proposition}
{\bf Proof.}
The proof is divided into 6 steps.
In Step 1, we analyze the geometric property of
the crossing of groups.
In Step 2, we make the assumptions about crossing.
In Step 3, we analyze the behavior of corigging
($=$ vacancy number $-$ rigging) under the operation of $\phi$.
Then we introduce a convenient graphical method to
analyze the coriggings
by using ``$\circ$" and ``$\bullet$".
Here ``$\bullet$" represents the letters 1 contained in
the local energy distribution and ``$\circ$"
represents the change of the quantum space
induced by letters 1 contained in the path.
In Step 4, we derive the relations from the assumption
made in Step 2.
In Step 5, we consider how to minimize the number of
``$\circ$" for given pattern of ``$\bullet$" of
the local energy distribution.
Finally, in Step 6, we show that combination of
the relations derived in Step 4 and
arguments in Step 5 lead to the contradiction,
hence completes the proof.
\vspace{3mm}

\noindent
{\it Step 1.}
Consider the path $b=b_1\otimes b_2\otimes\cdots\otimes b_L$,
and calculate the local energy distribution corresponding to
$b_{[k]}=b_1\otimes b_2\otimes\cdots\otimes b_k$ ($k\leq L$).
Recall that in Lemma \ref{s:lem:e_distrib->conf},
we have shown that patterns of the local energy distribution
represent the combinatorial procedures of the KKR
bijection $\phi$.
Suppose that there are two groups of letters 1
whose cardinalities are $m_1$ and $m_2$, respectively,
below $b_{[k]}$.
We name these two groups as $M_1$ and $M_2$, respectively.
Here we take the top end of the group $M_1$ is located
to the left of that of the group $M_2$.

{}From geometric property of crossing,
we show that, by appropriate choice of $k$,
we can assume $m_1<m_2$
without crossing beneath $b_{[k]}$.
Assume that there is a crossing between $l$th row
and $(l+1)$th row,
whereas there is no crossing above it.
Denote the elements of $M_1$ and $M_2$ at $l$th row
by $(l,m_{1-})$ and $(l,m_{2-})$ where $m_{1-}<m_{2-}$,
and the elements of $m_1$ and $m_2$ at ($l+1$)th row
by $(l+1,m_{1+})$ and $(l+1,m_{2+})$ where $m_{2+}<m_{1+}$,
respectively.
{}From the procedure given in Step 2 of Lemma \ref{s:junbi},
we have $m_{2-}\leq m_{2+}$.
Now we choose $k$ such that $m_{2+}\leq k<m_{1+}$ is satisfied.
Since $k$ satisfies $m_{1-}<m_{2-}\leq m_{2+}\leq k$,
cardinality of the group $M_1$ beneath $b_{[k]}$ is $l$, whereas
that of the group $M_2$ is equal to or greater than $l+1$,
which gives the claim.
\vspace{3mm}

\noindent
{\it Step 2.}
We keep the notation like $m_{i+}$ etc., as before,
therefore we have $m_1<m_2$ beneath $b_{[m_{1+}-1]}$.
Again, we are assuming that the crossing of $M_1$ and $M_2$ occurs
beneath $b_{[m_{1+}]}$,
and also that there is no crossing beneath $b_{[m_{1+}-1]}$.
We denote the cardinalities of
the groups $M_1$ and $M_2$ under $b_{[m_{1+}-1]}$
by $m_1$ and $m_2$, respectively.
At the end of the proof, we will show that
the existence of the crossing leads to contradiction.
We consider the crossing of two groups,
since this is the fundamental situation.
The general case involving more than two crossings
follows from this fundamental case.
The situation here is depicted in the following diagram.
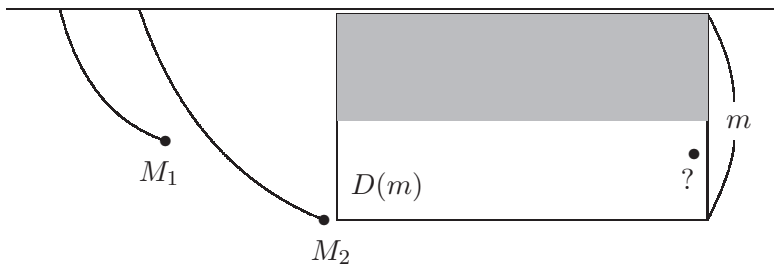
\begin{center}
\unitlength 10pt
\begin{picture}(29,11)
\put(0,10){\line(1,0){29}}
\qbezier(2,10)(3,6)(6,5)
\put(6,5){\circle*{0.4}}
\put(5,3.5){$M_1$}
\qbezier(5,10)(7,4)(12,2)
\put(12,2){\circle*{0.4}}
\put(11.5,0.5){$M_2$}
\multiput(12.5,2)(14,0){2}{\line(0,1){7.8}}
\multiput(12.5,2)(0,7.8){2}{\line(1,0){14}}
\put(13,3){$D(m)$}
\put(26,4.5){\circle*{0.4}}
\put(25.5,3.2){?}
\qbezier(26.5,9.8)(27.5,8)(27.5,6.5)
\qbezier(26.5,2)(27.5,4)(27.5,5)
\put(27.2,5.5){$m$}
\color[cmyk]{0,0,0,0.3}
\put(12.51,6){\rule[-2.2pt]{139.8pt}{40pt}}
\color{black}
\end{picture}
\end{center}
Here letters 1 in the local energy distribution
are represented by ``$\bullet$" and all letters 0
are suppressed.
Note that we have introduced the domain $D(m)$ on the right of
the bottom point of $M_2$, occupying from the first row
to the $m$th row.
Since we are considering the crossing caused by two groups,
we can assume that
groups contained in $D(m_2)$ are, in fact,
contained in $D(m_1-1)$ (in the above diagram,
it is indicated by the gray rectangle).
If we can put ``$\bullet$" indicated by ``?" in
the above diagram, then the crossing of the groups
$M_1$ and $M_2$ occurs.

In the following, we first treat the case that $M_2$ and
other groups on the right of $M_2$ are well separated.
This means the other groups on the right of $M_2$ are
located on the right of the bottom point of $M_2$,
and between the top and the bottom point of $M_2$,
there is no ``$\circ$" on the right of $M_2$
(see Step 3 for meaning of ``$\circ$").
This assumption is only for the sake of simplicity,
and the general case will be mentioned at the end of the proof.
\vspace{3mm}

\noindent
{\it Step 3.}
We summarize the basic properties of the vacancy numbers
(or, at the same time, that of the coriggings).
Recall the definition of the vacancy numbers
$p_j=Q^{(0)}_j-2Q^{(1)}_j$ corresponding to
the pair $(\lambda ,\mu)$.
Consider the box adding procedure of $\phi$.
If we add boxes to $\lambda |_{\leq j}$ and $\mu |_{\leq j}$
simultaneously, then the vacancy number
$p_j$ decrease by 1.
On the contrary, if we add a box to $\lambda |_{\leq j}$
and do not add box to $\mu |_{\leq j}$, then the vacancy
number increase by 1.
Note that If we do not add box to both $\lambda |_{\leq j}$
and $\mu |_{\leq j}$, then the vacancy number do not change.
Recall also that the procedure $\phi$ only refers
to corigging.

In order to analyze the above
change of coriggings induced by box
adding procedure of $\phi$, it is convenient to supplement
the local energy distribution with the information
of change of the quantum space corresponding to letters
1 contained in the path.
In the local energy distribution, we replace letters
1 by ``$\bullet$", and suppress all letters 0.
We join ``$\bullet$" belonging to the same group by thick lines.
Then, corresponding to the letters 1 contained in $b_{s}$,
we put ``$\circ$" on the right of ``$\bullet$"
corresponding to the letters 2 contained in $b_s$.
The row coordinates of ``$\circ$" are taken
to be the same as the column coordinate of the added
box of the quantum space induced by the corresponding letters 1.
Here we give examples of such diagram for
$\fbox{22}\otimes\fbox{12}\otimes\fbox{1}\otimes\fbox{22}$ and
$\fbox{22}\otimes\fbox{112}\otimes\fbox{22}$,
respectively.
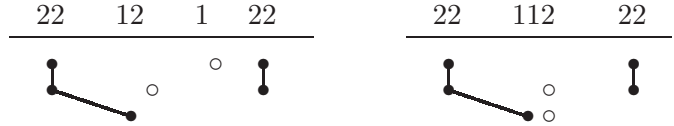
\begin{center}
\unitlength 10pt
\begin{picture}(25,5)
\put(0,3.5){\line(1,0){11.5}}
\put(1,4){22}
\put(4,4){12}
\put(7,4){1}
\put(9,4){22}
\put(1.6,2.5){\circle*{0.4}}
\put(1.6,1.5){\circle*{0.4}}
\put(5.4,1.5){\circle{0.4}}
\put(4.6,0.5){\circle*{0.4}}
\put(7.8,2.5){\circle{0.4}}
\put(9.6,2.5){\circle*{0.4}}
\put(9.6,1.5){\circle*{0.4}}
\thicklines
\put(1.6,1.5){\line(0,1){1}}
\qbezier(1.6,1.5)(1.6,1.5)(4.6,0.5)
\put(9.6,1.5){\line(0,1){1}}
\thinlines
\put(15,0){
\put(0,3.5){\line(1,0){10}}
\put(1,4){22}
\put(4,4){112}
\put(8,4){22}
\put(1.6,2.5){\circle*{0.4}}
\put(1.6,1.5){\circle*{0.4}}
\put(5.4,1.5){\circle{0.4}}
\put(5.4,0.5){\circle{0.4}}
\put(4.6,0.5){\circle*{0.4}}
\put(8.6,2.5){\circle*{0.4}}
\put(8.6,1.5){\circle*{0.4}}
\thicklines
\put(1.6,1.5){\line(0,1){1}}
\qbezier(1.6,1.5)(1.6,1.5)(4.6,0.5)
\put(8.6,1.5){\line(0,1){1}}
}
\end{picture}
\end{center}
In this diagram, the KKR map $\phi$ proceeds from the left to right,
and within the same column
(distinguishing columns of ``$\bullet$" and ``$\circ$"),
it proceeds from the top to bottom.
In the above examples, the left group containing three ``$\bullet$"
became singular (i.e., corigging $=0$)
after the bottom ``$\bullet$" is added.
Then the two ``$\circ$" increase the corigging by 2,
thereby the right group containing two ``$\bullet$"
stays independently from the left group.
This kind of analysis of change of the coriggings
is a prototype of the arguments given in Step 5 and Step 6.
Note that along each group of ``$\bullet$",
the notion of left and right of the group is well defined.
Let us remark
the convenient method to determine locations of ``$\circ$".
Given an element \fbox{$1\cdots 12\cdots 2$},
we reverse the orderings of numbers as
\fbox{$2\cdots 21\cdots 1$}.
Choose the specific letter 1, and denote by
$p$ the number of letters 1 and 2
on the left of it.
Then, corresponding to the chosen 1,
we put ``$\circ$" on the ($p+1$)th row on
the local energy distribution.
\vspace{3mm}

\noindent
{\it Step 4.}
Assume that we are going to add the box corresponding
to $(m_1+1,m_{1+})$ of the local energy distribution.
In order to add a box corresponding to 
$(m_1+1,m_{1+})$, or in other words, in order to
make crossing, the row of $\mu$ corresponding to $M_2$
cannot be singular when we add the box
corresponding to $(m_1+1,m_{1+})$.
This follows from the assumption $m_1<m_2$
and the fact that
we add a box to the longest possible singular row
in the procedure $\phi$.
Also, the row of $\mu$ corresponding to $M_1$ have to be singular
in order to add a box corresponding to $(m_1+1,m_{1+})$.
We consider the implications of these two conditions.

Recall that the row of $\mu$ corresponding to the group $M_2$
is singular when the bottom ``$\bullet$" is added to
the end of $M_2$.
On the other hand, we have to make $M_2$ non-singular
as we have seen in the above.
This means
\begin{equation}
(\mbox{number of ``$\circ$" within }D(m_2))>
(\mbox{number of ``$\bullet$" within }D(m_2)).
\end{equation}
On the other hand, in order to make $M_1$ singular, we have
\begin{equation}\label{s:circ_leq_bullet}
(\mbox{number of ``$\circ$" within }D(m_1))\leq
(\mbox{number of ``$\bullet$" within }D(m_1)).
\end{equation}
Note that all ``$\bullet$" are contained in $D(m_1-1)$.
{}From these two restrictions, we see there are at least one
``$\circ$" in $D(m_2)\setminus D(m_1)$.
\vspace{3mm}

\noindent
{\it Step 5.}
If we are given the pattern of ``$\bullet$",
there remains various possibilities about pattern of ``$\circ$".
Now we are going to consider the patterns of ``$\circ$"
that minimize the number of ``$\circ$".
To say the result at first, we see that we need
``$\circ$" as much as ``$\bullet$", therefore
in order to meet the condition (\ref{s:circ_leq_bullet}),
we have to minimize the number of ``$\circ$" .

Suppose there are two groups in the local energy distribution,
the group $s_1$ on the left, and $s_2$ on the right.
Let the top ``$\bullet$" be located at columns $k_1$ and $k_2$.
Then, in order to make $s_1$ and $s_2$ as separated groups,
we need $\min (s_1,s_2)$ ``$\circ$" within the region
between (or surrounded by) $s_1$ and column $k_2$.
To make the situation transparent, we consider the concrete path
$\fbox{222222222}\otimes\fbox{111111}\otimes\fbox{2}\otimes
\fbox{1122}\otimes\fbox{22}\otimes\fbox{22}\otimes
\fbox{111111}\otimes\fbox{2222222}$.
Then the corresponding diagram is as follows.
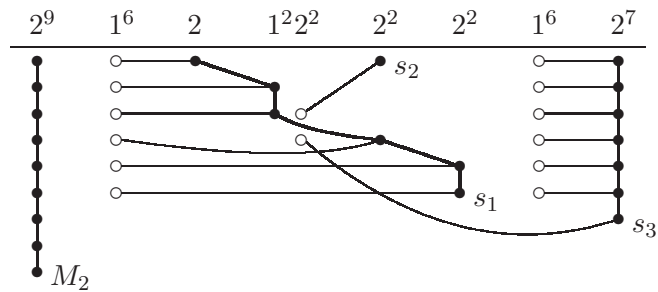
\begin{center}
\unitlength 10pt
\begin{picture}(24,11)
\put(-1,8.5){\line(1,0){24}}
\multiput(0,0)(0,1){9}{\circle*{0.4}}
\put(-0.3,9){$2^9$}
\multiput(3,3)(0,1){6}{\circle{0.4}}
\put(2.7,9){$1^6$}
\put(6,8){\circle*{0.4}}
\put(5.7,9){$2$}
\multiput(9,6)(0,1){2}{\circle*{0.4}}
\multiput(10,5)(0,1){2}{\circle{0.4}}
\put(8.7,9){$1^22^2$}
\multiput(13,5)(0,3){2}{\circle*{0.4}}
\put(12.7,9){$2^2$}
\multiput(16,3)(0,1){2}{\circle*{0.4}}
\put(15.7,9){$2^2$}
\multiput(19,3)(0,1){6}{\circle{0.4}}
\put(18.7,9){$1^6$}
\multiput(22,2)(0,1){7}{\circle*{0.4}}
\put(21.7,9){$2^7$}
\thicklines
\put(0,0){\line(0,1){8}}
\qbezier(6,8)(6,8)(9,7)
\put(9,6){\line(0,1){1}}
\qbezier(9,6)(10,5.3)(13,5)
\qbezier(13,5)(13,5)(16,4)
\put(16,3){\line(0,1){1}}
\put(22,2){\line(0,1){6}}
\thinlines
\put(3.2,8){\line(1,0){2.8}}
\put(3.2,7){\line(1,0){5.8}}
\put(3.2,6){\line(1,0){5.8}}
\qbezier(3.2,5)(10.2,4)(13,5)
\put(3.2,4){\line(1,0){12.8}}
\put(3.2,3){\line(1,0){12.8}}
\qbezier(10.2,6.1)(10.2,6.1)(13,8)
\multiput(19.2,3)(0,1){6}{\line(1,0){2.8}}
\qbezier(10.2,4.9)(16,0)(22,2)
\put(0.5,-0.5){$M_2$}
\put(16.5,2.5){$s_1$}
\put(13.5,7.5){$s_2$}
\put(22.5,1.5){$s_3$}
\end{picture}
\end{center}
We see there are 4 groups, labeled by $M_2$,
$s_1$, $s_2$, $s_3$ from left to right.
These groups are indicated by thick lines.
We can analyze the situation as follows
(as for the method for
analysis of change of the coriggings,
see latter part of Step 3):
\begin{enumerate}
\item
In order to make $M_2$ and $s_1$ separated,
we need 6 ``$\circ$" between $M_2$
and $s_1$.
Precise meaning of ``between $M_2$ and $s_1$" etc.,
is given after the example.
In this example, they are supplied by 6 ``$\circ$"
coming from the left \fbox{111111}.
This situation is indicated by thin lines,
which join the corresponding ``$\circ$" and ``$\bullet$".
Of course, there is ambiguity in the way of joining
``$\circ$" and ``$\bullet$",
however this ambiguity brings no important effect,
hence we neglected.
For example, we can join top 5 ``$\circ$" coming
from the left \fbox{111111} and the bottom
``$\circ$" coming from \fbox{1122} with
6 ``$\bullet$" of $s_1$.
In such a case, the bottom ``$\circ$" coming
from the left \fbox{111111} should be connected
with the bottom ``$\bullet$" of $s_3$.
\item
In order to make $s_1$ and $s_2$ separated,
we need 1 ``$\circ$" between $s_1$ and $s_2$.
In this example, it is supplied by 1 ``$\circ$"
coming from \fbox{1122}.
\item
In order to make $s_2$ and $s_3$ separated,
we need 1 ``$\circ$" between $s_2$ and $s_3$.
In this example, it is supplied by the top ``$\circ$"
coming from the right \fbox{111111}.
\item
In order to make $s_1$ and $s_3$ separated,
we need 5 ``$\circ$" between $s_1$ and $s_3$.
In this example, it is supplied by the bottom 5 ``$\circ$"
coming from the right \fbox{111111}.
\item
In order to make $M_2$ and $s_3$ separated,
we need 1 ``$\circ$" between $M_2$ and $s_3$.
In this example, it is supplied by 1 ``$\circ$"
coming from \fbox{1122}.
\end{enumerate}
Let us remark that if we move one 1 of the second term
of the above path to the seventh term, i.e.,
$\fbox{222222222}\otimes\fbox{11111}\otimes\fbox{2}\otimes
\fbox{1122}\otimes\fbox{22}\otimes\fbox{22}\otimes
\fbox{1111111}\otimes\fbox{2222222}$
has exactly the same pattern of ``$\bullet$" as the above example.

{}From this example, we can infer the general case.
Assume there are groups $s_1$, $\cdots$, $s_{n}$
(from left to right)
on the right of $M_2$.
We define the region between $s_i$ and $s_j$ ($i<j$) as the
region surrounded by $s_i$, the first row and
the row containing the bottom point of $s_i$,
$s_j$ and the column containing the bottom point of $s_j$
(except $s_i$, $s_j$ and the column containing
the bottom point of $s_j$, see the following diagram).
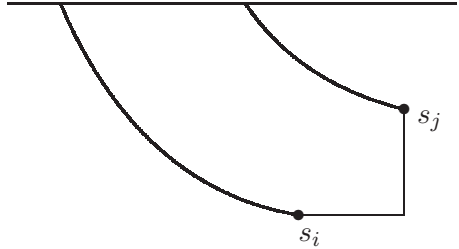
\begin{center}
\unitlength 10pt
\begin{picture}(15,10)(0,0.5)
\put(-1,10){\line(1,0){17}}
\thicklines
\qbezier(1,10)(4,3)(10,2)
\put(10,2){\circle*{0.4}}
\qbezier(8,10)(10,7)(14,6)
\put(14,6){\circle*{0.4}}
\thinlines
\put(10,2){\line(1,0){4}}
\put(14,2){\line(0,1){4}}
\put(10,1){$s_i$}
\put(14.5,5.5){$s_j$}
\end{picture}
\end{center}

We choose ordered subsequence
$s_{j_1}$, $\cdots$, $s_{j_p}$
of groups
$s_1$, $\cdots$, $s_{n}$
such that it is the longest subsequence which satisfies
$s_{j_1}>\cdots >s_{j_p}>s_{n-1}$.
Between $s_{n-1}$ and the $s_{n}$,
we need at least $\min (s_n,s_{n-1})$
``$\circ$"
in order to make $s_{n-1}$ and $s_{n}$ separate.
If $s_n>s_{n-1}$, then we compare $s_n$ and $s_{j_p}$.
Then we need at least $\min(s_n-s_{n-1},s_{j_p}-s_{n-1})$
``$\circ$" between $s_{j_p}$ and the $s_{n}$
in order to make $s_{j_p}$ and $s_n$ separate.
We continue this process and, in conclusion,
we need ``$\circ$" as much as ``$\bullet$",
compared within the right of $M_2$.
\vspace{3mm}

\noindent
{\it Step 6.}
Based on the ground of the arguments given in Step 5,
we derive the contradiction against the statement
``there are at least one ``$\circ$" in
$D(m_2)\setminus D(m_1)$".

Again assume there are groups $s_1$, $\cdots$, $s_{n}$
on the right of $M_2$.
Denote the column coordinate of the top ``$\bullet$"
of $s_j$ by $k_j$ and that of the bottom of
$s_j$ by $k_j'$.
In order to meet the condition (\ref{s:circ_leq_bullet}),
we choose the pattern of ``$\circ$" which minimizes the
number of ``$\circ$".
{}From observations made in Step 5,
we have at most $s_{n}$ ``$\circ$"
between $s_{n-1}$ and the $s_{n}$.
If all ``$\circ$" are located between the $k_{n-1}$th
column and $k_n$th column, then the number of
``$\circ$" (i.e. $s_{n}$) is
too short to make ``$\circ$" appear in
$D(m_2)\setminus D(m_1)$.
Let us analyze the case when some of $s_{n}$ ``$\circ$"
appear between the $k_n$th column and $k_n'$th column.
To be specific, take some $k$ between $k_n$ and $k_n'$,
and write the row of the lowest ``$\bullet$" of the column $k$
belonging to the group $s_n$ by $s$.
In order to make the upper $s$ ``$\bullet$" of $s_n$
separated from $s_{n-1}$, it needs at least $s$ ``$\circ$"
on the left of the $k$th column.
See the following schematic diagram.
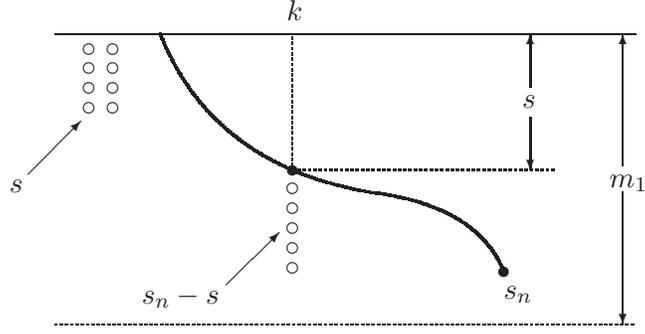
\begin{center}
\unitlength 10pt
\begin{picture}(22,12)(0,-1)
\put(0,10){\line(1,0){22}}
\thicklines
\qbezier(4,10)(6,5)(12,4)
\qbezier(12,4)(16,3.5)(17,1)
\put(17,1){\circle*{0.4}}
\put(17,0){$s_n$}
\thinlines
\multiput(9,5)(0,0.2){25}{\line(0,1){0.1}}
\multiput(9,4.85)(0.2,0){50}{\line(1,0){0.1}}
\put(8.8,10.5){$k$}
\put(9,4.85){\circle*{0.4}}
\multiput(1.3,7.2)(0,0.75){4}{\circle{0.4}}
\multiput(2.2,7.2)(0,0.75){4}{\circle{0.4}}
\put(-1,4.7){\vector(1,1){2}}
\put(-1.7,4.0){$s$}
\put(18,8){\vector(0,1){2}}
\put(18,6.85){\vector(0,-1){2}}
\put(17.7,7.2){$s$}
\multiput(9,4.15)(0,-0.75){5}{\circle{0.4}}
\put(6.5,0.5){\vector(1,1){2}}
\put(3.3,-0.2){$s_n-s$}
\multiput(0,-1)(0.2,0){110}{\line(1,0){0.1}}
\put(21.5,5){\vector(0,1){5}}
\put(21.5,4){\vector(0,-1){5}}
\put(21,4.2){$m_1$}
\end{picture}
\end{center}
Thus, we have at most $s_n-s$ ``$\circ$"
below the $s$th row of the $k$th column.
If we attach a column of $s_n-s$ ``$\circ$" to the $k$th row,
it has to begin from the $s+1$th row.
Therefore, it is also too short to make ``$\circ$" appear in
$D(m_2)\setminus D(m_1)$.
Similarly, between $s_{n-2}$ and the $s_{n-1}$,
we have at most $\max (s_{n},s_{n-1})$,
again too short to make ``$\circ$" appear in
$D(m_2)\setminus D(m_1)$.
Continuing in this way, we see that no ``$\circ$" appear in
$D(m_2)\setminus D(m_1)$, which gives contradiction.

As we have claimed at the end of Step 2, so far
we are dealing only with the case that $M_2$ and
the other groups on the right of it are well separated.
However, we can treat the general case by similar arguments.
First, by applying the same argument of Step 5,
we can show that we need ``$\circ$" as much as ``$\bullet$"
within the region on the right of the group $M_2$.
Then, by applying the same argument in Step 6,
all ``$\circ$" on the right of $M_2$ are included
in the first $m_1-1$ rows of the local energy distribution.
Therefore, in order to ``$\circ$" appear in
$D(m_2)\setminus D(m_1)$, we have to add at least one
``$\circ$" within the first $m_1$ rows.
This makes $M_1$ non-singular, hence the crossing
does not occur in this case.

Hence we complete the proof of proposition.
\hfill$\square$

\vspace{3mm}

\noindent
{\bf Proof of Theorem \ref{s:main}.}
From Proposition \ref{s:no-crossing}, we can remove
the procedure to find
``the nearest 1" from Step 2 of Lemma \ref{s:junbi}.
This gives the proof of Step 1 to Step 3 of Theorem \ref{s:main}.

Finally, we clarify the meaning of Step 4.
In Step 3 of Lemma \ref{s:junbi},
we determined $(\mu_k,j_k)$,
which corresponds to the rightmost box of row $\mu_k$ of $\mu$.
Since the row $\mu_k$ is not lengthened in calculation
of $\phi (b)$ after $b_{j_k}$,
the rigging of row $\mu_k$ is equal to the vacancy
number at the time when we add the rightmost box to $\mu_k$.
At this moment, the quantum space takes the form
\begin{equation}
(\lambda_1,\lambda_2,\cdots,\lambda_{j_k-1},E_{\mu_k,j_k}).
\end{equation}
The meaning of the last $E_{\mu_k,j_k}$ is as follows.
$E_{\mu_k,j_k}$ counts all letters 1 contained
in the first $\mu_k$ rows of $j_k$ th column of
local energy distribution.
This means that,
from Lemma \ref{s:lem:e_distrib->conf},
we added $E_{\mu_k,j_k}$ boxes
to $\mu$ before the rightmost box of the row $\mu_k$ is added
(while considering $b_{j_k}$).
In the procedure $\phi$, we use letters 2 of $b_{j_k}$ first,
and then use the rest of letters 1 of $b_{j_k}$.
Since letters 2 in $b_{j_k}$ means simultaneous
addition of box to the quantum space and $\mu$,
we can conclude that the quantum space
has the row with length $E_{\mu_k,j_k}(\leq\mu_k)$.
{}From this shape of the quantum space, we have
\begin{equation}
Q^{(0)}_{\mu_k}=
\sum_{i=1}^{j_k-1}\min (\mu_k,\lambda_i)
+E_{\mu_k,j_k}
\end{equation}

{}From Lemma \ref{s:lem:E=Q} applying to the path
$b_1\otimes b_2\otimes\cdots\otimes b_{j_k}$,
we deduce the following:
\begin{equation}
Q^{(1)}_{\mu_k}=
\sum_{i=1}^{j_k}E_{\mu_k,i}.
\end{equation}
Hence we obtain the formula in Step 4, and
complete the proof of Theorem \ref{s:main}.
\hfill$\square$

\vspace{3mm}

\noindent
{\bf Proof of Theorem \ref{s:thm:takagi}.}
This follows immediately from non-crossing property
of Proposition \ref{s:no-crossing}.
In Theorem \ref{s:main}, we determine groups from right to left.
More precisely,
after determining one group, all letters 1 belonging
to the group are made to be 0, and we determine
the rightmost group again.
However, if we start from the bottom point of the longest
group, the elimination made in the right has no effect,
and we can determine the same group by virtue of
non-crossing property.
By doing these procedure from longer groups to shorter groups,
we obtain the same groups as in the case Theorem \ref{s:main}.
\hfill$\square$

\section{Summary}\label{sec:summary}
In this paper, we consider crystal interpretation of
the KKR map $\phi$ from paths to rigged configurations.
In Section \ref{s:sec:main}, we introduce table
called the local energy distribution.
The entries of the table are differences of the energy
functions, and we show in Theorem \ref{s:main}
that this table have complete
information about $\phi$ so that we can read off the rigged
configuration from it.
As we see in Lemma \ref{s:lem:e_distrib->conf},
this table can be viewed as giving crystal interpretation
of the combinatorial procedures appearing in
the original definition of $\phi$.

As we see in Proposition \ref{s:no-crossing},
our formalism has simple property.
This enables us to reformulate Theorem \ref{s:main}
as described in Theorem \ref{s:thm:takagi}.
The latter formalism is particularly important when we
consider inverse scattering formalism for the
periodic box-ball systems.
As we see in Remark \ref{s:rem:pbbs},
advantage of our formalism,
compared with the formalism given in \cite{KTT}, is that we can treat
states of the periodic box-ball system directly
without sending them to linear semi-infinite systems.

\appendix
\section{Kerov--Kirillov--Reshetikhin bijection}\label{sec:kkr}
In this section, we prepare notations and basic properties
corresponding to the Kerov--Kirillov--Reshetikhin (KKR)
bijection.
As for the definitions of the rigged configurations
as well as
combinatorial procedure of the bijection, we refer
to Section 2 of \cite{Sak} ($\phi$ there should be read as
$\phi^{-1}$ here)
or Appendix A to \cite{KTT},
and we only prepare necessary notations.

Assume we have given a highest weight path $b$:
\begin{equation}
b=b_1\otimes b_2\otimes\cdots\otimes b_L\in
B_{\lambda_1}\otimes B_{\lambda_2}\otimes\cdots\otimes B_{\lambda_L}.
\end{equation}
Then we have one to one correspondence $\phi$ between $b$ and the
rigged configuration
\begin{equation}
\phi :b
\longrightarrow {\rm RC}=\left( (\lambda_i)_{i=1}^L,\,
(\mu_i,r_i)_{i=1}^N\right) .
\end{equation}
Here $(\mu_i)_i\in\mathbb{Z}_{\geq 0}^N$ is called
configuration and we depict $(\lambda_i)_i$ and $(\mu_i)_i$
by Young diagrammatic expression
whose rows are given by $\lambda_i$ and $\mu_i$, respectively.
Integers $r_i$ are called riggings and we associate them with
the corresponding row $\mu_i$.
In the KKR bijection, orderings within integer sequences
$(\lambda_i)_i$
or $(\mu_i,r_i)_i$ does not make any differences.
On this rigged configurations, we use the symbols
$Q_j^{(a)}$ ($a=0,1$) defined by 
\begin{equation}\label{def:Q}
Q_j^{(0)}:=\sum_{k=1}^{L}\min (j,\lambda_k),
\qquad
Q_j^{(1)}:=\sum_{k=1}^{N}\min (j,\mu_k).
\end{equation}
The vacancy number $p_j$ for length $j$ rows of
$\mu$ is then defined by 
\begin{equation}
p_j:=Q^{(0)}_j-2Q^{(1)}_j.
\end{equation}
If row $\mu_i$ has property $p_{\mu_i}=r_i$,
then the row is called singular.
For the highest weight paths, the corresponding rigged
configurations are known to satisfy 
$0\leq r_i\leq p_{\mu_i}$.
The quantity $p_{\mu_i}-r_i$ is sometimes called corigging.

One of the most important properties of $\phi$ or
$\phi^{-1}$ is that if we consider isomorphic paths
$b\simeq b'$, then the corresponding rigged configuration
is the same (Lemma 8.5 of \cite{KSS}).
We express this property in terms of the map $\phi^{-1}$
as follows.
\begin{theorem}\label{s:thm:KSS}
Take successive two rows from the quantum space $\lambda$
of the rigged configuration arbitrary,
and denote them by $\lambda_a$ and $\lambda_b$.
When we remove $\lambda_a$ at first
and next $\lambda_b$ by the KKR map $\phi^{-1}$,
then we obtain two tableaux,
which we denote by $a_1$ and $b_1$, respectively.
Next, on the contrary, we first remove $\lambda_b$
and second $\lambda_a$ (keeping the order of other removal invariant)
and we get $b_2$ and $a_2$.
Then we have
\begin{equation}
b_1\otimes a_1\,\simeq\, a_2\otimes b_2 ,
\end{equation}
under the isomorphism of $\mathfrak{sl}_2$
combinatorial $R$ matrix.
\hfill$\square$
\end{theorem}

We remark that there is an extension of $\phi$ which covers
all non-highest weight elements as well.
Let $b$ be an arbitrary element of arbitrary
tensor products of crystals;
$b\in B_{\lambda_1}\otimes\cdots\otimes B_{\lambda_L}$.
In particular, $b$ can be non-highest weight element.
Then we can apply the same combinatorial procedure
for $\phi$ and obtain $\phi (b)$
as extension of the rigged configurations.
Following \cite{Sch,DS}, we call
such $\phi (b)$ {\it unrestricted rigged configuration}.
Let us denote $\phi(b)=(\lambda,(\mu,r))$.
Then, from definition of $\phi$, we see that
$|\lambda |$ represents the number of all letters 1 and
2 contained in the path $b$, whereas $|\mu |$
represents the number of letters 2 contained in $b$.
Note, in particular, that $|\lambda |\geq |\mu |$
holds for unrestricted rigged configurations.
These unrestricted rigged configurations contain
the rigged configurations as the special case.

Let $b$ be a non-highest weight element as above.
Consider the following modification of $b$:
\begin{equation}
b':=\fbox{1}^{\,\otimes \Lambda}\otimes b,
\end{equation}
where $\Lambda$ is an integer satisfying
$\Lambda\geq \lambda_1+\cdots +\lambda_L$.
Then $b'$ is highest weight.
Under these notations, we have the following:

\begin{lemma}\label{s:lem:haiesutoka}
Let the unrestricted rigged configuration
corresponding to $b$ be
\begin{equation}\label{s:eq:unristrictedrc}
\bigl( (\lambda_i)_{i=1}^L ,(\mu_j ,r_j)_{j=1}^N\bigl).
\end{equation}
Then the rigged configuration corresponding to
the highest path $b'$ is given by
\begin{equation}\label{s:eq:unristrictedrc_modify}
\bigl( (\lambda_i)_{i=1}^L\cup (1^\Lambda)
,(\mu_j ,r_j+\Lambda)_{j=1}^N\bigl).
\end{equation}
\end{lemma}
{\bf Proof.}
Let the vacancy number of row $\mu_j$
of the pair $(\lambda ,\mu )$ of (\ref{s:eq:unristrictedrc})
be $p_{\mu_j}$.
Then the vacancy number of row $\mu_j$ of
(\ref{s:eq:unristrictedrc_modify})
is equal to $p_{\mu_j}+\Lambda$,
because of the addition of $(1^\Lambda)$ on $\lambda$.
Now we apply $\phi^{-1}$ on
(\ref{s:eq:unristrictedrc_modify}).
{}From $\lambda\cup (1^\Lambda)$ of the quantum space,
we remove $\lambda$ first, and next remove $(1^\Lambda)$.
Recall that in the combinatorial procedure of $\phi^{-1}$,
we only refer to corigging,
and it does not refer to actual values of the riggings.
Therefore, when we remove $\lambda$ from the quantum space
of (\ref{s:eq:unristrictedrc_modify}), we obtain
$b$ as the corresponding part of the image.
Then, remaining rigged configuration has
the quantum space $(1^\Lambda)$ without $\mu$ part.
On this rigged configuration,
the map $\phi^{-1}$ becomes trivial
and obtain $b'$ as the image corresponding to
(\ref{s:eq:unristrictedrc_modify}).
\hfill$\square$

\section{Operators $T_l$}\label{sec:T_l}
In this section, we introduce the operators $T_l$
which are defined by the combinatorial $R$.
$T_l$'s serve as the time evolution operators
of the box-ball systems \cite{HHIKTT}.
Let $u_l$ be a highest weight element of $B_l$,
i.e., in a tableau representation, it is
$u_l=\fbox{$\overbrace{11\cdots 1}^{l}$}$.
We consider the path
\begin{equation}\label{s:eq:shapeofpath}
b=b_1\otimes b_2\otimes\cdots\otimes b_L\in
B_{\lambda_1}\otimes B_{\lambda_2}\otimes\cdots\otimes B_{\lambda_L}.
\end{equation}
Then its time evolution $T_l(b)$ $(l\in\mathbb{Z}_{>0})$
is defined by successively sending $u_l$ to the
right of $b$ under the isomorphism of the combinatorial
$R$ as follows:
\begin{eqnarray}
u_l\otimes b&=&
u_l\otimes b_1\otimes b_2\otimes\cdots\otimes b_L
\nonumber\\
&\stackrel{R}{\simeq}&
b_1'\otimes u_l^{(1)}\otimes b_2\otimes\cdots\otimes b_L
\nonumber\\
&\stackrel{R}{\simeq}&
b_1'\otimes b_2'\otimes u_l^{(2)}\otimes\cdots\otimes b_L
\nonumber\\
&\stackrel{R}{\simeq}&\cdots\cdots
\nonumber\\
&\stackrel{R}{\simeq}&
b_1'\otimes b_2'\otimes\cdots\otimes b_L'\otimes u_l^{(L)}
\nonumber\\
&=:&T_l(b)\otimes u_l^{(L)}.\label{s:def:T_l}
\end{eqnarray}

According to Proposition 2.6 of \cite{KOSTY}, operators $T_l$
on highest paths can be linearized by the
KKR bijection.
Since, in the main text, we use the similar property
for general case including non-highest paths,
we include here a proof for generalized version.

\begin{proposition}\label{s:prop:ist}
(1) Consider the path $b$ of the form (\ref{s:eq:shapeofpath}).
Here $b$ can be non-highest weight element.
Define $b'=b\otimes\fbox{1}^{\,\otimes\Lambda}$,
where the integer $\Lambda$ satisfies
$\Lambda >\lambda_1+\lambda_2+\cdots +\lambda_L$.
Then, we have $u_l\otimes b'\simeq T_l(b')\otimes u_l$.

(2) Denote the (unrestricted) rigged configuration corresponding
to $b'$ as
\begin{equation}\label{s:prop:T_l:1}
b'\xrightarrow{\mathrm{KKR}}
\bigl( \lambda\cup (1^\Lambda)
,(\mu_j ,r_j)_{j=1}^N\bigl).
\end{equation}
Then, corresponding to $T_l(b')$, we have
\begin{equation}\label{s:prop:T_l:2}
T_l(b')\xrightarrow{\mathrm{KKR}}
\bigl( \lambda\cup (1^\Lambda)
,(\mu_j ,r_j+\min (\mu_j,l))_{j=1}^N\bigl).
\end{equation}
\end{proposition}
{\bf Proof.}
(1) Consider the elements $u_l^{(i)}$ defined in
(\ref{s:def:T_l}).
In our case, $u_l^{(L)}$ contains letters 2 for at most
$\lambda_1+\lambda_2+\cdots +\lambda_L$ times.
Then, by calculating combinatorial $R$ along (\ref{s:def:T_l})
with $u_l^{(L)}$ and $\fbox{1}^{\,\otimes\Lambda}$,
we see that $u_l^{(L+\Lambda )}=u_l$.

(2) Consider the following rigged configuration:
\begin{equation}\label{s:prop:T_l:3}
\bigl( \lambda\cup (1^\Lambda)\cup (l)
,(\mu_j ,r_j+\min (\mu_j,l))_{j=1}^N\bigl),
\end{equation}
i.e., we added a row with width $l$ to the quantum space.
Compare the coriggings
of (\ref{s:prop:T_l:3}) and (\ref{s:prop:T_l:1}).
Recall that the vacancy number is defined by
$Q_{\mu_j}^{(0)}-2Q_{\mu_j}^{(1)}$.
As for $Q_{\mu_j}^{(1)}$, both (\ref{s:prop:T_l:3}) and
(\ref{s:prop:T_l:1}) give the same value,
since we have $\mu$ in the both second terms.
On the contrary, $Q_{\mu_j}^{(0)}$ for (\ref{s:prop:T_l:3})
is greater than the one for (\ref{s:prop:T_l:1})
by $\min (\mu_j,l)$, since we have the extra row of width
$l$ in the quantum space of (\ref{s:prop:T_l:3}).
In (\ref{s:prop:T_l:3}), riggings are increased by
value $\min (\mu_j,l)$, therefore we conclude that
the coriggings for both (\ref{s:prop:T_l:3})
and (\ref{s:prop:T_l:1}) coincide.

Now we apply $\phi^{-1}$ on (\ref{s:prop:T_l:3})
in two different ways.
First, we remove $\lambda\cup (1^\Lambda)$
from the quantum space of (\ref{s:prop:T_l:3})
(order of removal is the same as $\phi^{-1}$ on
(\ref{s:prop:T_l:1}) to obtain $b'$).
Since the coriggings for both $\mu_j$ of (\ref{s:prop:T_l:3})
and (\ref{s:prop:T_l:1}) coincide, we obtain
$b'$ as the corresponding part of the image.
Then we are left with the rigged configuration
$(l,(\emptyset ,\emptyset))$, which yields $u_l$.
Therefore we obtain $u_l\otimes b'$ as the image.

Next, we apply $\phi^{-1}$ on (\ref{s:prop:T_l:3})
in different way.
In this case, we remove the row $l$ of the quantum space
as the first step.
Note that in the (unrestricted) rigged configuration
$\bigl( \lambda ,(\mu ,r)\bigl)$
corresponding to the path $b$, all riggings $r_j$ are
smaller than the corresponding vacancy numbers.
By definition of $\Lambda$, we have
$\Lambda >\min (\mu_j,l)$ for all $j$
(recall that from definition of the unrestricted
rigged configuration, we always have
$\lambda_1+\cdots +\lambda_L\geq\mu_1+\cdots +\mu_N$).
As the result, if we remove the row $l$ from the
quantum space of (\ref{s:prop:T_l:3}),
rows $\mu_j$ do not become singular even if the riggings
are increased as $r_j+\min (\mu_j,l)$,
since the vacancy numbers are also increased by $\Lambda$
by the addition of $(1^\Lambda)$.
Thus, we obtain $u_l$ as the corresponding part of the image.
Then we are left with
$\bigl( \lambda\cup (1^\Lambda)
,(\mu_j ,r_j+\min (\mu_j,l))_{j=1}^N\bigl)$,
whose corresponding path we denote by $\tilde{b'}$.
In conclusion, we obtain $\tilde{b'}\otimes u_l$ as the image.

In the above two calculation of $\phi^{-1}$, the only
difference is the order of removing rows of the quantum
space of (\ref{s:prop:T_l:3}).
Therefore we can apply Theorem \ref{s:thm:KSS}
to get the isomorphism
\begin{equation}
u_l\otimes b'\stackrel{R}{\simeq}
\tilde{b'}\otimes u_l.
\end{equation}
If $b'$ is non-highest weight, we apply Lemma
\ref{s:lem:haiesutoka} and
$u_l\otimes u_m\simeq u_m\otimes u_l$,
then we can use the same argument.
{}From (1), this means $\tilde{b'}=T_l(b')$,
which completes a proof.
\hfill$\square$

\vspace{6mm}

\noindent
{\bf Acknowledgements:}
The author would like to thank Atsuo Kuniba
and Taichiro Takagi for collaboration in an
early stage of the present study
and valuable comments on this manuscript.
He is a research fellow of the 
Japan Society for the Promotion of Science.

\end{document}